\documentclass[12pt,leqno]{article}
\usepackage{amsthm,amsfonts,amssymb,epsfig,graphics,amsmath,eufrak}
 \usepackage[latin1]{inputenc}\relax

\newcommand{\R}{\Re}

\newcommand\adots{\mathinner{\mkern2mu\raise1pt\hbox{.}
\mkern3mu\raise4pt\hbox{.}\mkern1mu\raise7pt\hbox{.}}}

\newcommand{\ubar}{\bar{u}}
\newcommand{\utild}{\tilde{u}}
\newcommand{\fe}{\varphi}
\newcommand{\hcal}{\mathcal{H}}
\newcommand{\fcal}{\mathcal{F}}
\newcommand{\gtild}{\tilde{G}}

\newcommand{\dcal}{\mathcal{D}}

\newcommand{\ddp}{\frac{\partial\ubar^\delta}{\partial
\delta}}

\newcommand{\mathbbR}{\mathbb{R}}
\newcommand{\ta}{\tilde{A}}

\newtheorem{theo}{Theorem}[section]
\newtheorem{prop}[theo]{Proposition}
\newtheorem{cor}[theo]{Corollary}
\newtheorem{lem}[theo]{Lemma}
\newtheorem{defi}[theo]{Definition}

\newtheorem{rem}[theo]{Remark}

\numberwithin{equation}{section}

\title{Pointwise Asymptotic Behavior of Perturbed Viscous Shock Profiles}

\author{\sc \small
Peter Howard\thanks{Texas A\& M University, College Station, TX
77843; phoward@math.tamu.edu} \, and Mohammadreza
Raoofi\thanks{Indiana University, Bloomington, IN 47405;
mraoofi@indiana.edu. } }
\begin{document}

\maketitle
\begin{abstract}
We consider the asymptotic behavior of perturbations of Lax and
overcompressive type viscous shock profiles arising in systems
of regularized conservation laws with strictly parabolic viscosity,
and also in systems of conservation laws with partially parabolic
regularizations such as arise in the case of the compressible
Navier--Stokes equations and in the equations of magnetohydrodynamics.
Under the necessary conditions of spectral and hyperbolic stability,
together with transversality of the connecting profile, we establish
detailed pointwise estimates on perturbations from a sum of the viscous
shock profile under consideration and a family of diffusion waves
which propagate perturbation signals
along outgoing characteristics.  Our approach combines the recent
$L^p$-space analysis of Raoofi \cite{Ra} with a straightforward
bootstrapping argument that relies on a refined
description of nonlinear signal interactions, which we develop through
convolution estimates involving Green's functions for the linear
evolutionary PDE that arises upon linearization of the
regularized conservation law about the distinguished profile.
Our estimates are similar to, though slightly weaker than,
those developed by Liu in his landmark result on the case of weak Lax type
profiles arising in the case of identity viscosity \cite{Liu97}.
\end{abstract}

\bigbreak
\section{Introduction}

\medskip

Consider a ``viscous shock profile'', or traveling-wave solution
\begin{equation}
 u(x,t)=\bar u(x-st), \qquad
 \lim_{z\to \pm \infty} \bar u(z)= u_\pm,
\label{profile}
\end{equation}
of a second-order system of conservation laws
\begin{equation}
\begin{aligned}
 u_t+ F(u)_{x}&=
(B(u)u_{x})_{x}, \cr x\in \mathbb{R}; \,  &u, \, F\in \mathbb{R}^n;
\, B\in  \mathbb{R}^{n\times n}, \label{viscous}
\end{aligned}
\end{equation}
for which we take
\begin{equation}
u=\left(\begin{matrix}  u^I\\ u^{II}\end{matrix}\right), \quad
B=\left(\begin{matrix} 0 & 0\\
b_1 & b_2 \end{matrix}\right), \label{UB}
\end{equation}
$ u^I\in \mathbbR^{n-r}$, $u^{II}\in \mathbbR^r$,
and
\begin{equation}
Re \sigma  b_2 \ge \theta , \label{goodb}
\end{equation}
with $\theta>0$ and where $\sigma$ denotes spectrum.

Equations of form (\ref{viscous}) arise in a wide variety of
applications, including such well known examples as the compressible
Navier--Stokes equations, the equations of magnetohydrodynamics (see
\cite{Z.5}), and the equations of viscoelastic materials.  The long
time behavior of solutions of (\ref{viscous}) is often characterized
by a family of viscous shock profiles (\ref{profile}), which serves
as a regularization for the shock pattern solution of the associated
hyperbolic problem
\begin{equation}
 u_t + F(u)_x = 0.
 \label{hyperbolic}
\end{equation}
In this context, we expect that a
viscous profile will persist in the pattern only if it
is individually stable to small fluctuations in initial
conditions, and consequently the stability of viscous shock
profiles has long been a subject of considerable interest and
research effort (see \cite{ZH} and the references therein).

An important advance in the study of such waves and their
stability was the identification in \cite{ZH} of a spectral
criterion based on the {\it Evans function} (see particularly
\cite{AGJ, E, GZ, J, KS, MaZ.3, SZ, Z.4, ZH}).  Briefly,
the Evans function, typically denoted $D(\lambda)$,
serves as a characteristic
function for the linear operator $L$ that arises upon
linearization of (\ref{viscous}) about the stationary
profile $\bar{u} (x)$.  More
precisely, away from essential spectrum, zeros of the
Evans function correspond in location and multiplicity with
eigenvalues of $L$ \cite{AGJ, GZ, ZH}.  It was shown in
\cite{ZH} and \cite{MaZ.3}, respectively for the
strictly parabolic and real viscosity cases, that
$L^1 \cap L^p \to L^p$ linearized orbital stability of the
profile, $p>1$, is equivalent to the Evans function
condition,
\medbreak
($\mathcal{D}$) \quad There exist
precisely $\ell$ zeroes of $D(\cdot)$ in the nonstable half-plane
$\R \lambda \ge 0$, necessarily at the origin $\lambda=0$.
\medbreak
\noindent
Here, $\ell$ is the dimension of the manifold connecting
$u_-$ and $u_+$.

Stability criterion $(\mathcal{D})$ has been shown to hold in all
cases for small amplitude Lax shocks arising in strictly parabolic
systems \cite{FreS, Go.1, HuZ, KM, KMN, MN, PZ}, as well as for
large-amplitude shocks in such cases as Lax type waves arising in
isentropic Navier--Stokes equations for the gamma-law gas as $\gamma
\to 1$ \cite{MN}, and undercompressive shocks arising in Slemrod's
model for van der Waal gas dynamics \cite{Z.5} (see \cite{S} for
Slemrod's model). More generally, condition $(\mathcal{D})$ can be
verified by numerical calculation \cite{B.1, B.2, BZ, BDG}. In the
case of strictly parabolic systems such {\it spectral stability},
along with standard technical hypotheses on $F$ and $B$ (see
$(\hcal0)$--$(\hcal4)$, in Remark \ref{strictpc} below,)  has been
shown sufficient for establishing nonlinear stability of Lax,
under-, over-, and mixed under--overcompressive shock profiles
\cite{HZ}, while for mixed hyperbolic--parabolic regularization,
condition $\mathcal{D}$ (along with (H0)--(H3) and (A1)--(A3) below)
has been shown sufficient for establishing nonlinear stability for
Lax and overcompressive shock profiles \cite{Ra, MaZ.4, Z.5}. Except
in the case of \cite{Ra}, these nonlinear analyses are carried out
through consideration of the perturbation
\begin{equation*}
v(x,t) = u(x,t) - \bar{u}^{\delta(t)} (x),
\end{equation*}
where $\delta(t)$ is introduced as a local phase chosen to insure that
the shapes of $u$ and $\bar{u}$ are compared, rather than their positions
(in this way, orbital stability is considered).  In this context,
and under the assumption of stability criterion $(\mathcal{D})$,
it can be shown that for initial perturbations
\begin{equation*}
|v(x,0)|\le E_0 (1+|x|)^{-3/2},
\end{equation*}
some $E_0$ sufficiently small, there holds
\begin{equation} \label{introLp}
\|v(x,t)\|_{L^p} \le C E_0 (1+t)^{-\frac{1}{2} (1-\frac{1}{p})},
\end{equation}
from which we observe asymptotic decay in time for all $p > 1$.

A natural refinement of this type of analysis regards the
consideration of {\it diffusion waves}, which are defined as
exact solutions to a family of convecting Burgers' equations
(convection along outgoing characteristics of the underlying
hyperbolic problem),
and which carry precisely the $L^1$ mass in (\ref{introLp})
that does not decay asymptotically in time
\cite{Liu85, Liu97}.  Considering, then, the perturbation
\begin{equation*}
v(x,t) = u(x,t) - \bar{u}^{\delta(t)} (x) - \varphi (x,t),
\end{equation*}
where $\varphi (t,x)$ represents the sum of diffusion waves,
we have
\begin{equation*}
\int_{-\infty}^{+\infty} \varphi (x,0) dx =
\int_{-\infty}^{+\infty} (u(x,0) - \bar{u}^{\delta(t)} (x)) dx,
\end{equation*}
and consequently,
\begin{equation*}
\int_{-\infty}^{+\infty} v (x,0) dx = 0.
\end{equation*}
In this way, we reduce to the case of zero-mass initial data, for
which perturbations decay equally from above and below the profile,
and the asymptotic rate of decay is doubled (at least at the linear
level, and for an initial perturbation with sufficiently fast
spatial decay).  Working in this context, Raoofi has shown that
stability criterion $(\mathcal{D})$, along with (H0)--(H3) and
(A1)--(A3), are sufficient for establishing that for initial
perturbations
\begin{equation*}
v(x,t) = u(x,t) - \bar{u}^{\delta_*} (x) - \varphi (x,0) -
\frac{\partial \bar{u}^\delta}{\partial \delta} \delta (t),
\end{equation*}
where $\delta_*$ is the asymptotic shape/location of the shock, and
with $|u(\cdot,0) - \bar{u}|$ sufficiently small, there holds
\cite{Ra}
\begin{equation*}
\|v (\cdot,t)\|_{L^p}
\le C E_0 (1+t)^{-\frac{1}{2}(1-\frac{1}{p})-\frac{1}{4}}.
\end{equation*}

The goal of the current analysis is both to refine the
analysis of Raoofi to the case of pointwise (rather than
$L^p$) estimates, and to introduce a simplified
bootstrapping argument through which estimates on
the perturbation $v(t,x)$ emerge in straightforward
fashion.  Our estimates are similar to those of
Liu, developed for Lax type profiles in the case of identity
viscosity and under the assumption of weak shock
strength.  Our analysis has no such limitations, though
we note that the form of our viscosity and the shock
strength under consideration are encoded in our
spectral assumptions.

Throughout the analysis, we will work in a coordinate
system moving along with the shock, so that without loss
of generality, we consider a standing profile $\ubar (x)$,
which satisfies
the traveling-wave ordinary differential equation
(ODE)
\begin{equation}
B(\bar u) \bar u'=F(\bar u)-F(u_-). \label{ODE}
\end{equation}
Considering the block structure of $B$, this can be written as:
\begin{equation}
F^I(u^I, u^{II})\equiv F^I(u_-^I, u_-^{II})\label{eq1}
\end{equation}
and
\begin{equation}
b_1(u^I)' + b_2(u^{II})'= F^{II}(u^I, u^{II}) - F^{II}(u_-^I,
u_-^{II}). \label{eq2}
\end{equation}
We are interested in the asymptotic and pointwise behavior of
$\utild$, a solution of (\ref{viscous}) and a perturbation of
$\ubar$.

\medskip
 We assume that, by some invertible change of coordinates
$u\to w(u)$, possibly but not necessarily connected with a global
convex entropy, followed if necessary by multiplication on the left
by a nonsingular matrix function $S(w)$, equations (\ref{viscous})
may be written in the {\it quasilinear, partially symmetric
hyperbolic-parabolic form}
\begin{equation}
\tilde A^0 w_t +  \tilde A w_{x}= (\tilde B w_{x})_{x} + G, \quad
w=\left(\begin{matrix} w^I
\\w^{II}\end{matrix}\right), \label{symm}
\end{equation}
 $w^I\in \mathbbR^{n-r}$, $w^{II}\in \mathbbR^r$, $x\in
\mathbbR^d$, $t\in \mathbbR$, where, defining $W_\pm:= W(U_\pm)$:\\
\medskip
(A1)\quad $\tilde A(w_\pm)$, $\tilde A_{11}$, $\tilde
A^0$ are symmetric, $\tilde A^0 >0$.\\
\medskip
(A2)\quad Dissipativity: no eigenvector of $ dF(u_\pm)$ lies in the
kernel of $ B(u_\pm)$. (Equivalently, no eigenvector of $ \tilde A
(\tilde A^0)^{-1} (w_\pm)$
lies in the kernel of $\tilde B (\tilde A^0)^{-1}(w_\pm)$.)\\
\medskip
(A3) \quad $ \tilde B= \left(\begin{matrix} 0 & 0 \\ 0 & \tilde b
\end{matrix}\right) $, $ \tilde G= \left(\begin{matrix}  0 \\ \tilde g\end{matrix}\right) $,
with $ Re  \tilde b(w) \ge \theta $ for some $\theta>0$, for all
$W$, and $\tilde g(w_x,w_x)=\mathbf{O}(|w_x|^2)$.
\medskip
Here, the coefficients of (\ref{symm}) may be expressed in terms of
the original equation (\ref{viscous}), the coordinate change $u\to
w(u)$, and the approximate symmetrizer $S(w)$, as
\begin{equation}
\begin{aligned} \tilde A^0&:= S(w)(\partial u/\partial w),\quad
\tilde A:= S(w)d F(u(w))(\partial u/\partial w),\\
\tilde B&:= S(w)B(u(w))(\partial u/\partial w), \quad G= -(dS w_{x})
B(u(w))(\partial u/\partial w) w_{x}.
\end{aligned}
\label{coeffs}
\end{equation}

\smallskip
Along with the above structural assumptions, we make
the technical
hypotheses:\\
\medskip
(H0)\quad $F$, $B$, $w$, $S\in C^{8}$.\\
\smallskip
(H1)\quad The eigenvalues of $\tilde A_*:=\ta_{11}(\ta^0_{11})^{-1}$
are (i) distinct from $0$; (ii) of common sign; and (iii) of
constant multiplicity with respect to
$u$.\\
\smallskip
(H2)\quad The eigenvalues of $dF(u_\pm)$ are real, distinct, and
nonzero.\\
\smallskip
({H3})\quad Nearby $\bar u$, the set of all solutions of
\eqref{profile}--\eqref{viscous} connecting the same values $u_\pm$
forms a smooth manifold $\{\bar u^\delta\}$, $\delta\in
\mathcal{U}\subset \mathbb{R}^\ell$, $\bar u^0=\bar u$.

We note that structural assumptions (A1)--(A3) and technical
hypotheses (H0)--(H2) are broad enough to include such cases
as the compressible Navier--Stokes equations, the equations
of magnetohydrodymics, and Slemrod's model for van der Waal gas
dynamics \cite{Z.5}.   Moreover, existence of waves $\bar{u}$
satisfying (H3) has been established in each of these cases.

\bigskip
\begin{defi}\label{type}
An ideal shock
\begin{equation} \label{shock}
u(x,t) =
\begin{cases}
u_- &x < st, \\
u_+ &x > st,
\end{cases}
\end{equation}
is classified as
{\it undercompressive}, {\it Lax},
or {\it overcompressive} type according as $i-n$
is less than, equal to, or greater than $1$,
where $i$, denoting the sum of the dimensions $i_-$ and $i_+$
of the center--unstable subspace of $df(u_-)$ and the center--stable
subspace of $df(u_+)$, represents the total number of characteristics
incoming to the shock.

A viscous profile \eqref{profile} is classified
as {\it pure undercompressive} type if the associated ideal shock
is undercompressive and $\ell=1$,
{\it pure Lax} type
if the corresponding ideal shock is Lax type and $\ell=i-n$,
and {\it pure overcompressive} type if
the corresponding ideal shock is overcompressive and $\ell=i-n$,
$\ell$ as in (H3).
Otherwise it is classified as {\it mixed under--overcompressive} type;
see \cite{LZ.2, ZH}.
\end{defi}

Pure Lax type profiles are the most common type,
and the only type
arising in standard gas dynamics, while pure over- and undercompressive type
profiles arise in magnetohydrodynamics (MHD) and phase-transitional models.

Under assumptions (A0)--(A3) and (H0)--(H3), or their analogs in the
real viscosity case, condition ($\mathcal{D}$) is equivalent to (i)
{\it strong spectral stability}, $\sigma(L)\subset \{\R \lambda \le
0\}\cup \{0\}$, (ii) {\it hyperbolic stability} of the associated
ideal shock, and (iii) {\it transversality} of $\bar u$ as a
solution of the connection problem in the associated traveling-wave
ODE, where hyperbolic stability is defined for Lax and
undercompressive shocks by the Lopatinski condition of \cite{M.1,
M.2, M.3, Fre} and for overcompressive shocks by an analogous
long-wave stability condition (see ($\mathcal{D}$ii) below); see
\cite{ZH, MaZ.3, ZS, Z.3, Z.4, Z.5} for further explanation. From
now on, we assume ($\mathcal{D}$) to hold along with (A0)--(A3) and
(H0)--(H4). We also assume that the the shock is a pure Lax or
overcompressive one.

\medskip
Setting $A_\pm:=df(u_\pm)$, $\Gamma_\pm:=d^2f(u_\pm)$, and $B_\pm:=B(u_\pm)$,
denote by
\begin{equation}
a_1^-< a_2^-<\dots < a_n^- \quad \text{\rm and }
a_1^+< a_2^+<\dots < a_n^+
\label{a}
\end{equation}
the eigenvalues of $A_-$ and $A_+$, and $l_j^\pm$, $r_j^\pm$ left
and right eigenvectors associated with each $a_j^\pm$, normalized so
that $(l_j^{T}r_k)_\pm=\delta^j_k$, where $\delta^j_k$ is the
Kronecker delta function, returning $1$ for $j=k$ and $0$ for $j\ne
k$. Under this notation,  hyperbolic stability of $\ubar$, a Lax or
overcompressive shock profile,  is the condition:

\medbreak ($\mathcal{D}$ii)\quad The set $\{r^\pm; a^\pm \gtrless
0\} \cup \{\int_{-\infty}^{+\infty}\frac{\partial
\ubar^\delta}{\partial \delta_i} dx; i=1, \cdots, \ell \}$ forms a
basis for $\mathbb{R}^n$, with
$\int_{-\infty}^{+\infty}\frac{\partial \ubar^\delta}{\partial
\delta_i} dx$  computed at $\delta=0.$ \medskip \\ As said before,
($\mathcal{D}$ii) is satisfied whenever ($\mathcal{D}$) holds.

 Define scalar diffusion coefficients
\begin{equation}
\beta_j^\pm:= (l_j^T Br_j)_\pm
\label{beta}
\end{equation}
and scalar coupling coefficients
\begin{equation}
\gamma_j^\pm:= (l_j^T \Gamma (r_j,r_j))_\pm.
\label{gamma}
\end{equation}

Following \cite{Liu, Liu85, Liu97},
define for a given mass $m_j^-$ the
scalar diffusion waves $\fe_j^-(x,t;m_j^-)$ as
(self-similar) solutions of the Burgers equations
\begin{equation}
\fe_{j,t}^- + a^-_j\fe^-_{j,x} - \beta^-_j \fe_{j,xx}^- = -\gamma^-_j
((\fe_{j}^-)^2)_x
\label{burgers}
\end{equation}
with point-source initial data
\begin{equation}
\fe_j^- (x, -1) = m_j^-\delta_0(x),
\label{burgersdata}
\end{equation}
and similarly for $\fe_j^+(x,t;m_j^+)$.
Given a collection of masses
$m_j^\pm$
prescribed on outgoing characteristic modes
$a_j^-<0$ and $a_j^+>0$, define
\begin{equation}
\fe(x,t) =\sum_{a_j^- <0}\fe_j^-(x,t; m_j^-) r^-_j
+ \sum_{a_j^+>0}\fe_j^+(x,t; m_j^+) r^+_j.
\label{phi}
\end{equation}

Also define

\begin{equation}\label{psi1}
\begin{aligned}
\psi_1(x,t)&:=\chi(x,t) \sum_{a_j^-<0}
(1+t)^{-1/2} (1+|x-a_j^-t|+t^{\frac13})^{-3/4}\\
&\quad+ \chi(x,t)\sum_{a_j^+>0}
(1+t)^{-1/2} (1+|x-a_j^+t|+t^{\frac13})^{-3/4},\\
\end{aligned}
\end{equation}
and
\begin{equation}\label{psi1bar}
\begin{aligned}
\bar\psi_1(x,t)&:=\chi(x,t) \sum_{a_j^-<0}
(1+t)^{-1/2} (1+|x-a_j^-t|)^{-3/4}\\
&\quad+ \chi(x,t)\sum_{a_j^+>0}
(1+t)^{-1/2} (1+|x-a_j^+t|)^{-3/4},\\
\end{aligned}
\end{equation}
where $\chi(x,t)=1$ for $x\in [a_1^-t, a_n^+t]$ and zero otherwise.
Also,
\begin{equation}\label{psi2}
\begin{aligned}
\psi_2(x,t)&:=  \sum_{a_j^-<0}
 (1+|x-a_j^-t|+t^{1/2})^{-3/2}\\
&\quad+ \sum_{a_j^+>0} (1+|x-a_j^+t|+t^{1/2})^{-3/2},\\
\end{aligned}
\end{equation}
and
\begin{equation}
\alpha(x,t) := \chi (x,t) (1+t)^{-3/4}(1+|x|)^{-1/2}.
\end{equation}

  Theorem  \ref{thmmain} and  Corollaries \ref{cormain} and \ref{corliu} are the main results of this paper:
\begin{theo} \label{thmmain}
 Assume (A1)--(A3), (H0)--(H3) and $(\mathcal{D})$ hold, and  $\ubar$ is a pure Lax or overcompressive shock profile.
  Assume also that $\utild$ solves (\ref{viscous}) with
initial data $\utild_0$ and that, for initial perturbation
$u_0:=\utild_0-\ubar$, we have $|u_0|_{L^1\cap H^4} \leq E_0$,
$|u_0(x)|\le E_0(1+|x|)^{-\frac32},$ and $|\partial_x u_0(x)|\le
E_0(1+|x|)^{-\frac12},$ for $E_0$ sufficiently small. Then there are
an $\ell$-array  function $\delta(t)$, and a small constant
$\delta_*$, such that if
$v:=\utild-\ubar^{\delta_*}-\fe-\ddp\delta$, then
\begin{equation}
|v(x,t)|\le
 C E_0
(\psi_1+\psi_2 +\alpha), \label{1stestimate}
\end{equation}
and
\begin{equation}
 |v_x(x,t)|\le
 C E_0
t^{-\frac12}(1+t)^{\frac12}(\bar\psi_1+\psi_2 +\alpha);
\label{2stestimate}
\end{equation}
furthermore,
\begin{equation} \label{deltabound}
|\delta (t)|\le C E_0(1+t)^{-\frac 12},
\end{equation}
and
\begin{equation}\label{dotdeltabd}
|\dot \delta (t)|\le C E_0(1+t)^{-1},
\end{equation}
 for some constant $C$ (independent of $x,t$ and $E_0$).
\end{theo}

The proof of Theorem \ref{thmmain}
uses a straightforward bootstrapping argument, combined
with $L^p$ estimates proved in \cite{Ra} and restated here in the
following proposition.

\begin{prop}\label{propra} Under the conditions of Theorem \ref{thmmain},
\begin{equation}\label{vlp}
|v(\cdot,t)|_{L^p}\le
 C E_0
(1+t)^{-\frac{1}{2}(1-\frac1p)-\frac 14},
\end{equation}
for any $p,$ $1 \leq p \leq \infty$; also
\begin{equation}\label{vhs}
|v(\cdot,t)|_{H^s}\le
 C E_0
(1+t)^{-\frac{1}{2}},
\end{equation}
for any $s\le 4$, and for some constant $C$.
\end{prop}
\noindent {\bf Proof of Proposition \ref{propra}.} See \cite{Ra}.
In \cite{Ra} the bound (\ref{vhs}) was established  for $s\le 3$;
basically, to prove (\ref{vhs}), it was shown, using energy
estimates, that the $H^3$-norm of $u$ is controlled by the
$H^3$-norm of the initial data and the $L^2$-norm of $u$. However,
the same proof (with a slight modification) can be used for $s=4;$
we need only to assume that the initial data is small in $H^4$ and
the coefficients are in $C^8$ (rather than $C^6$ as in \cite{Ra}).
The reason we need two times differentiability of the coefficients
is in the fact that we need the exponential decay of $\ubar$ and
its derivatives to their endstates up to $2s$ derivatives; see
(5.56) in \cite{Ra} and Lemma \ref{expdecay} in the present work.
\hfill $\square$

\bigskip  Taylor expansion gives us
\begin{equation}
\ubar^{\delta_*+\delta(t)}-\ubar^{\delta_*}=\ddp\delta(t) +
\mathbf O (|\delta(t)|^2 e^{-k|x|}),\label{ishl}
\end{equation}
\medskip
since $\ubar$ approaches in an exponential rate to its endstates
(see Lemma \ref{expdecay}). But $|\delta(t)|^2 e^{-k|x|}$ is then
smaller than the right hand side of (\ref{1stestimate}). Hence we
have the following.
\begin{cor} \label{cormain} Under the assumption of  Theorem
\ref{thmmain},
\begin{equation}
\utild(x,t)-\ubar^{\delta_*+\delta(t)}(x)=\fe(x,t)+
\mathbf{O}(\psi_1+\psi_2 +\alpha).
\end{equation}
Also,
\begin{equation}
|\utild(\cdot,t)-\ubar^{\delta_*+\delta(t)}-\fe(\cdot,t)|_{L^p}=
\mathbf{O}((1+t)^{-\frac12(1-\frac1p)-\frac14}).
\end{equation}
\end{cor}

Without ``instantaneous shock tracking" $\delta(t)$, however, we
obtain the following, which Liu  proved for the \emph{artificial
viscosity} case \cite{Liu97}.
\begin{cor} \label{corliu} Under the assumption of  Theorem
\ref{thmmain},
\begin{equation}
|\utild(\cdot,t)-\ubar^{\delta_*}-\fe(\cdot,t)|_{L^p}=\begin{cases}
\mathbf{O}((1+t)^{-\frac12(1-\frac1p)-\frac14})\quad &\textup{for}\, 1\le p\le 2\\
\mathbf{O}((1+t)^{-\frac12})\quad &\textup{for}\, 2\le p\le \infty.
\end{cases}
\end{equation}
\end{cor}

\medskip
The picture of asymptotic behavior described in Theorem
\ref{thmmain} and corollary \ref{cormain} was introduced on
heuristic grounds by Liu \cite{Liu85} in the context of
small-amplitude Lax-type shock waves and artificial (identity)
viscosity $B=I$, and, along with the accompanying analysis of
\cite{Liu} described below, played an important role in the
subsequent analysis by Szepessy and Xin in \cite{SX} establishing
for the first time stability (with no rate) of small-amplitude
Lax-type shock profiles with $B=I$. Bounds \eqref{vlp} validating
this picture were established for large-amplitude Lax or
overcompressive profiles and general viscosity by Raoofi \cite{Ra}.
See also earlier pointwise arguments sketched (but not completed) in
\cite{Liu97,ZH}. In the proof of Theorem \ref{thmmain} we will
assume and use the bounds (\ref{vlp}) and (\ref{vhs}), already
established in \cite{Ra}.

\begin{rem} \textup{
The estimates of Theorem \ref{thmmain} are similar to, though
weaker than, those developed by Liu in the case of
identity viscosity and weak shock strength \cite{Liu97}.
In particular, Liu's estimates provide slightly sharper results in terms
of the pointwise bounds, and also give more information along the
characteristic modes. In order to provide results as sharp as those
of Liu, we would have to treat the critical nonlinearity
$\varphi v$ similarly as we treat the nonlinearity $\varphi^2$,
which we analyze by a clever approach of Liu's in which he
integrates the nonlinear interaction integral by parts in
$t$, making use of the observation that
$(\partial_t + a_j^- \partial_x) \varphi_j^-$ decays
similarly as the $t$-derivative of a heat kernel \cite{Liu}.
(For the argument in our context,
see (\ref{diffmain}) and the surrounding discussion, for which
under a shift of coordinates, differentiation with respect
to $\tau$ replaces the operator $(\partial_t + a_j^- \partial_x)$.
The resulting estimate is stated in (\ref{gphi}) of
Lemma \ref{linear}.)
In order to apply this approach to the term $\varphi v$,
we would additionally have to
carry an estimate on $(\partial_t + a_j^- \partial_x) v$
through our argument.  We leave the full details of this
calculation to a separate paper \cite{HRZ}.
We remark also that the estimates of Liu are uniform
in shock strength $\epsilon := |u_+ - u_-|$ as
$\epsilon \to 0$.  More specifically, $\epsilon$ appears
explicitly in Liu's estimates, and the estimates remain
largely unchanged as $\epsilon \to 0$, though one
of the estimates (denoted $\chi_i$) increases slightly
due to the loss of a term $\epsilon^{-1} \psi_i (x, t)$
as a possibility in an estimate that takes the minimum of
three quantities (the other two quantities remain uniform
in $\epsilon$) \cite{Liu97}.  In our case, shock strength is assumed
fixed, and our estimates are not uniform as $\epsilon \to 0$.
In particular, coefficients in the Green's function estimates
we use blow up as $\epsilon \to 0$, and we must counter this
by reducing $E_0$, the small constant multiplying our initial
perturbation.  We regard the extension of our analysis to
the case uniform in (small) shock strength as an
important future project.
One of the advantages
of our approach, on the other hand, is the use of the instantaneous
shock tracking $\delta(t)$, as well as using the previously
established $L^p$ norms, which makes our proof both more
straightforward and easier to generalize.  We also remark
that since we proceed from Green's function estimates
obtained for possibly large-amplitude shock profiles
(information about the amplitude is encoded in the spectrum
of the linearized operator), our analysis applies to this
case as well.
Finally, the $L^p$ estimates of Corollary \ref{corliu} recover
the estimates of \cite{Ra}, with a slight improvement
on the shift location estimates $|\delta (t)|$ and $\dot \delta (t)$.
As mentioned in \cite{Ra}, these estimates, in the
case $p > 2$, are a slight improvement of Liu's.
}
\end{rem}
\begin{rem}[Remarks on the Strictly Parabolic case] \label{strictpc}
\textup{  The case of the strictly parabolic systems can be treated
in a very similar way, the proof being almost identical. However, we
need fewer assumptions for the equations, or for the initial
perturbation, in this case. Basically, instead of (A1)--(A3) and
(H0)--(H3), we assume the following assumptions in the case of
strictly parabolic case.}
\medskip

({$\hcal0$}) \quad  $f,B\in C^3$.
\medskip

({$\hcal1$}) \quad $Re \, \sigma(B) > 0$.
\medskip

({$\hcal2$}) \quad  \textup{The eigenvalues of $df(u_\pm)$ are real,
distinct, and nonzero.}
\medskip

({$\hcal3$}) \quad  \textup{For some $\theta>0$, and all real $k$,
we have
$$Re\, \sigma(-ik df(u_\pm) -k^2 B(u_\pm))< -\theta k^2.$$ }

({$\hcal4$}) \quad \textup{The set of all stationary solutions of
(\ref{profile})-(\ref{viscous}) near $\ubar$, connecting the same
values $u_\pm$, forms a smooth  manifold $\{\ubar^{\delta}\},
\delta\in \mathbb{R}^{\ell}, \ubar^0=\ubar.$ }

\textup{We also assume that $(\mathcal{D})$ holds. For the initial
perturbation we only need $u_0 \in C^{1+\alpha}$ for some
$0<\alpha <1$.  In Theorem \ref{thmmain} concerning the real
viscosity (partially parabolic) case, we require the pointwise
bound on the derivative of the initial perturbation---i.e. the
bound $|\partial_x u_0(x)|\le E_0(1+|x|)^{-\frac12}$---only in
order to control the derivative in the hyperbolic modes, which are
absent in the strictly parabolic case. The necessary bounds on the
derivatives in the strictly parabolic case can be achieved through
short time estimates similar to ones carried out in \cite{Ra, HZ,
ZH}. On the other hand, as we are assuming less on the initial data,
(\ref{vhs}) does not necessarily hold. See \cite{Ra} for more
details.}

\end{rem}

\section{Linearized equations and Green function bounds}
Before stating the Green function bounds for the linearized
equation, we need some preparation to do. First we need the
exponential decay of $\ubar$ to its endstates. The following lemma
proved in \cite{MaZ.2} provides us with that.
\begin{lem} \label{expdecay}
Given (H1)--(H3), the endstates $u_\pm$ are hyperbolic rest points
of the ODE determined by (\ref{eq2}) on the $r$-dimensional manifold
(\ref{eq1}), i.e., the coefficients of the linearized equations
about $u_\pm$, written in local coordinates, have no center
subspace. In particular, under regularity (H0),
\begin{equation}
D_x^j D_\delta^i(\bar
u^\delta(x)-u_{\pm})=\mathbf{O}(e^{-\alpha|x|}), \quad \alpha>0, \,
0\le j\le 8, \,i=0,1,
\end{equation}
as $x\rightarrow\pm\infty$.
\end{lem}

%
%

Instead of linearizing about $\ubar(\cdot)$, we linearize about
$\ubar^{\delta_*}(\cdot)$, where $\delta_*$, determined \emph{a
priori} by the mass of the perturbation, would be the asymptotic
location or shape of the shock. Linearizing around
$\ubar^{\delta_*}(\cdot)$ gives us
\begin{equation}
v_t=Lv:=-(Av)_x+(Bv_x)_x, \label{linearov}
\end{equation}
with
\begin{equation}
B(x):= B(\ubar^{\delta_*}(x)), \quad  A(x)v:=
df(\ubar^{\delta_*}(x))v-dB(\ubar^{\delta_*}(x))v\ubar^{\delta_*}_x.
\label{AandBov}
\end{equation}
and let
\begin{equation}
G(x,t;y):= e^{Lt}\delta_y (x). \label{ov3.7}
\end{equation}

Denoting $A^\pm := A(\pm \infty)$,  $B^\pm:= B(\pm \infty)$, and
considering  lemma \ref{expdecay}, it follows that
\begin{equation}
|A(x)-A^-|= \mathbf {O} (e^{-\eta |x|}), \quad |B(x)-B^-|= \mathbf
{O} (e^{-\eta |x|}) \label{ABboundsov}
\end{equation}
as $x\to -\infty,$ for some positive $\eta.$ Similarly for $A^+$ and
$B^+,$ as $x\to +\infty.$ Also $|A(x)-A^\pm|$ and $|B(x)-B^\pm|$ are
bounded for all $x$.

 Define the {\it (scalar) characteristic speeds} $a^\pm_1<
\cdots < a_n^\pm$ (as above) to be the eigenvalues of $A^\pm$, and
the left and right {\it (scalar) characteristic modes} $l_j^\pm$,
$r_j^\pm$ to be corresponding left and right eigenvectors,
respectively (i.e., $A^\pm r_j^\pm = a_j^\pm r_j^\pm,$ etc.),
normalized so that $l^+_j \cdot r^+_k=\delta^j_k$ and $l^-_j \cdot
r^-_k=\delta^j_k$. Following Kawashima \cite{Kaw}, define associated
{\it effective scalar diffusion rates} $\beta^\pm_j:j=1,\cdots,n$ by
relation
\begin{equation}
\left(
\begin{matrix}
\beta_1^\pm &&0\\
&\vdots &\\
0&&\beta_n^\pm
\end{matrix}
\right) \quad = \hbox{diag}\ L^\pm B^\pm R^\pm, \label{betaov}
\end{equation}
where $L^\pm:=(l_1^\pm,\dots,l_n^\pm)^t$,
$R^\pm:=(r_1^\pm,\dots,r_n^\pm)$ diagonalize $A^\pm$.

Assume for  $A$ and $B$ the block structures:
$$A=\left(\begin{matrix}A_{11}\quad A_{12}\\A_{21}\quad A_{22}\end{matrix}\right),
B=\left(\begin{matrix}0& 0\\B_{21}& B_{22}\end{matrix}\right).$$


Also, let $a^{*}_j(x)$, $j=1,\dots,(n-r)$ denote the eigenvalues of
$$
A_{*}:= A_{11}- A_{12} B_{22}^{-1}B_{21}, \label{A*}
$$
with $l^*_j(x)$, $r^*_j(x)\in \mathbbR^{n-r}$ associated left and
right eigenvectors, normalized so that $l^{*t}_jr_j\equiv
\delta^j_k$. More generally, for an $m_j^*$-fold eigenvalue, we
choose $(n-r)\times m_j^* $ blocks $L_j^*$ and $R_j^*$ of
eigenvectors satisfying the {\it dynamical normalization}
$$
L_j^{*t}\partial_x R_j^{*}\equiv 0,
$$
along with the usual static normalization $L^{*t}_jR_j\equiv
\delta^j_kI_{m_j^*}$; as shown in Lemma 4.9, \cite{MaZ.1}, this may
always be achieved with bounded $L_j^*$, $R_j^*$. Associated with
$L_j^*$, $R_j^*$, define extended, $n\times m_j^*$ blocks
$$
\mathcal{L}_j^*:=\left(\begin{matrix} L_j^* \\
0\end{matrix}\right), \quad \mathcal{R}_j^*:=
\left(\begin{matrix} R_j^*\\
-B_{22}^{-1}B_{21} R_j^*\end{matrix}\right). \label{CalLR}
$$
%
Eigenvalues $a_j^*$ and eigenmodes $\mathcal{L}_j^*$,
$\mathcal{R}_j^*$ correspond, respectively, to short-time hyperbolic
characteristic speeds and modes of propagation for the reduced,
hyperbolic part of degenerate system (\ref{viscous}).

Define local, $m_j\times m_j$ {\it dissipation coefficients}
$$
\eta_j^*(x):= -L_j^{*t} D_* R_j^* (x), \quad j=1,\dots,J\le n-r,
\label{eta}
$$
where
$$
\aligned
&{D_*}(x):= \\
& \, A_{12}B_{22}^{-1} \Big[A_{21}-A_{22} B_{22}^{-1} B_{21}+ A_{*}
B_{22}^{-1} B_{21} + B_{22}\partial_x (B_{22}^{-1} B_{21})\Big]
\endaligned
\label{D*}
$$
is an effective dissipation analogous to the effective diffusion
predicted by formal, Chapman--Enskog expansion in the (dual)
relaxation case.

At $x=\pm \infty$, these reduce to the corresponding quantities
identified by Zeng \cite{Ze.1,LZe} in her study by Fourier transform
techniques of decay to {\it constant solutions} $\bar U \equiv
u_\pm$ of hyperbolic--parabolic systems, i.e., of limiting equations
$$
U_t=L_\pm U:= -A_{\pm} U_x+ B_\pm U_{xx}. \label{limiting}
$$
As a consequence of dissipativity, (A2), we obtain (see, e.g.,
\cite{Kaw, LZe, MaZ.3})
\begin{equation}
\beta_j^{\pm}>0, \quad Re \sigma(\eta_j^{*\pm})>0 \quad \text{\rm
for all $j$}. \label{goodbeta}
\end{equation}
However, note that the dynamical dissipation coefficient $D_*(x)$
{\it does not} agree with its static counterpart, possessing an
additional term $B_{22}\partial_x (B_{22}^{-1} B_{21})$, and so we
cannot conclude that (\ref{goodbeta}) holds everywhere along the
profile, but only at the endpoints. This is an important difference
in the variable-coefficient case; see Remarks 1.11-1.12 of
\cite{MaZ.3} for further discussion.

\begin{prop} \label{greenbounds}\cite{MaZ.3}  Under assumptions (A1)--(A3),
(H0)--(H3), and $(\mathcal D)$,
 the Green distribution $G(x,t;y)$
associated with the linearized evolution equations  may be
decomposed as
\begin{equation}
G(x,t;y)= H + E+  S + R, \label{ourdecomp}
\end{equation}
where, for $y\le 0$:
\begin{equation}
\begin{aligned} H(x,t;y)&:= \sum_{j=1}^{J} a_j^{*-1}(x) a_j^{*}(y)
\mathcal{R}_j^*(x) \zeta_j^*(y,t) \delta_{x-\bar a_j^* t}(-y)
\mathcal{L}_j^{*t}(y)\\
&= \sum_{j=1}^{J} \mathcal{R}_j^*(x) \mathcal{O}(e^{-\eta_0 t})
\delta_{x-\bar a_j^* t}(-y) \mathcal{L}_j^{*t}(y),
\end{aligned}
\label{multH}
\end{equation}
where the averaged convection rates $\bar a_j^*= \bar a_j^*(x,t)$ in
(\ref{multH}) denote the time-averages over $[0,t]$ of $a_j^*(x)$
along backward characteristic paths $z_j^*=z_j^*(x,t)$ defined by
$$
dz_j^*/dt= a_j^*(z_j^*), \quad z_j^*(t)=x, \label{char}
$$
and the dissipation matrix $\zeta_j^*=\zeta_j^*(x,t)\in
\mathbbR^{m_j^*\times m_j^*}$ is defined by the {\it dissipative
flow}
$$
d\zeta_j^*/dt= -\eta_j^*(z_j^*)\zeta_j^*, \quad
\zeta_j^*(0)=I_{m_j}; \label{diss}
$$
and $\delta_{x-\bar a_j^* t}$ denotes  Dirac distribution centered
at $x-\bar a_j^* t$.
\begin{equation}\label{E}
E(x,t;y):=\sum_{j=1}^\ell \frac{\partial \bar
u^\delta(x)}{\partial \delta_j}_{|\delta=\delta_*}e_j(y,t),
\end{equation}
\begin{equation}\label{e}
  e_j(y,t):=\sum_{a_k^{-}>0}
  \left(\textrm{errfn }\left(\frac{y+a_k^{-}t}{\sqrt{4\beta_k^{-}t}}\right)
  -\textrm{errfn }\left(\frac{y-a_k^{-}t}{\sqrt{4\beta_k^{-}t}}\right)\right)
  l_{jk}^{-}
\end{equation}
\begin{equation}
\begin{aligned} S(x,t;y)&:= \chi_{\{t\ge 1\}}\sum_{a_k^-<0}r_k^-  {l_k^-}^t
(4\pi \beta_k^-t)^{-1/2} e^{-(x-y-a_k^-t)^2 / 4\beta_k^-t} \\
&+\chi_{\{t\ge 1\}} \sum_{a_k^- > 0} r_k^-  {l_k^-}^t (4\pi
\beta_k^-t)^{-1/2} e^{-(x-y-a_k^-t)^2 / 4\beta_k^-t}
\left({\frac {e^{-x}}{e^x+e^{-x}}}\right)\\
&+ \chi_{\{t\ge 1\}}\sum_{a_k^- > 0, \,  a_j^- < 0}
[c^{j,-}_{k,-}]r_j^-  {l_k^-}^t (4\pi \bar\beta_{jk}^- t)^{-1/2}
e^{-(x-z_{jk}^-)^2 / 4\bar\beta_{jk}^- t}
\left({\frac{e^{ -x}}{e^x+e^{-x}}}\right),\\
&+\chi_{\{t\ge 1\}} \sum_{a_k^- > 0, \,  a_j^+ > 0}
[c^{j,+}_{k,-}]r_j^+  {l_k^-}^t (4\pi \bar\beta_{jk}^+ t)^{-1/2}
e^{-(x-z_{jk}^+)^2 / 4\bar\beta_{jk}^+ t}
\left({\frac{e^{ x}}{e^x+e^{-x}}}\right),\\
\end{aligned}
\label{Sov}
\end{equation}
with
\begin{equation}
z_{jk}^\pm(y,t):=a_j^\pm\left(t-\frac{|y|}{|a_k^-|}\right)
\label{zjkov}
\end{equation}
and
\begin{equation}
\bar \beta^\pm_{jk}(x,t;y):= \frac{x^\pm}{a_j^\pm t} \beta_j^\pm +
\frac{|y|}{|a_k^- t|} \left( \frac{a_j^\pm}{a_k^-}\right)^2
\beta_k^-, \label{betaaverageov}
\end{equation}
The remainder $R$ and its derivatives have the following bounds.
\begin{equation}
\begin{aligned}
R(x,t;y)&=
\mathbf{O}(e^{-\eta(|x-y|+t)})\\
&+\sum_{k=1}^n \mathbf{O} \left( (t+1)^{-1/2} e^{-\eta x^+}
+e^{-\eta|x|} \right)
t^{-1/2}e^{-(x-y-a_k^{-} t)^2/Mt} \\
&+ \sum_{a_k^{-} > 0, \, a_j^{-} < 0} \chi_{\{ |a_k^{-} t|\ge |y|
\}} \mathbf{O} ((t+1)^{-1/2} t^{-1/2})
e^{-(x-a_j^{-}(t-|y/a_k^{-}|))^2/Mt}
e^{-\eta x^+}, \\
&+ \sum_{a_k^{-} > 0, \, a_j^{+}> 0} \chi_{\{ |a_k^{-} t|\ge |y| \}}
\mathbf{O} ((t+1)^{-1/2} t^{-1/2}) e^{-(x-a_j^{+}
(t-|y/a_k^{-}|))^2/Mt}
e^{-\eta x^-}, \\
\end{aligned}
\label{Rbounds}
\end{equation}
\begin{equation}
\begin{aligned}
R_y(x,t;y)&= \sum_{j=1}^J \mathbf{O}(e^{-\eta t})\delta_{x-\bar
a_j^* t}(-y) +
\mathbf{O}(e^{-\eta(|x-y|+t)})\\
&+\sum_{k=1}^n \mathbf{O} \left( (t+1)^{-1/2} e^{-\eta x^+}
+e^{-\eta|x|}+e^{-\eta|y|}  \right)
t^{-1}
e^{-(x-y-a_k^{-} t)^2/Mt} \\
&+ \sum_{a_k^{-} > 0, \, a_j^{-} < 0} \chi_{\{ |a_k^{-} t|\ge |y|
\}} \mathbf{O} ((t+1)^{-1/2} t^{-1})
e^{-(x-a_j^{-}(t-|y/a_k^{-}|))^2/Mt}
e^{-\eta x^+} \\
&+ \sum_{a_k^{-} > 0, \, a_j^{+} > 0} \chi_{\{ |a_k^{-} t|\ge |y|
\}} \mathbf{O} ((t+1)^{-1/2} t^{-1})
e^{-(x-a_j^{+}(t-|y/a_k^{-}|))^2/Mt}
e^{-\eta x^-}, \\
\end{aligned}
\label{Rybounds}
\end{equation}
\begin{equation}
\begin{aligned}
R_x(x,t;y)&= \sum_{j=1}^J \mathbf{O}(e^{-\eta t})\delta_{x-\bar
a_j^* t}(-y) +
\mathbf{O}(e^{-\eta(|x-y|+t)})\\
&+\sum_{k=1}^n \mathbf {O} \left( (t+1)^{-1} e^{-\eta x^+}
+e^{-\eta|x|} \right)
t^{-1} (t+1)^{1/2}
e^{-(x-y-a_k^{-} t)^2/Mt} \\
&+ \sum_{a_k^{-} > 0, \, a_j^{-} < 0} \chi_{\{ |a_k^{-} t|\ge |y|
\}} \mathbf{O}(t+1)^{-1/2} t^{-1})
e^{-(x-a_j^{-}(t-|y/a_k^-|))^2/Mt}
e^{-\eta x^+} \\
&+ \sum_{a_k^{-} > 0, \, a_j^{+} > 0} \chi_{\{ |a_k^{-} t|\ge |y|
\}} \mathbf{O}(t+1)^{-1/2} t^{-1})
e^{-(x-a_j^{+}(t-|y/a_k^{-}|))^2/Mt}
e^{-\eta x^-}. \\
\end{aligned}
\label{Rxbounds}
\end{equation}
Moreover, for $|x-y|/t$ sufficiently large, $|G|\le Ce^{-\eta
t}e^{-|x-y|^2/Mt)}$ as in the strictly parabolic case.
\end{prop}
Setting $\gtild := S+R$, so that $G= H + E + \gtild$,  we have the
following alternative bounds for $\gtild$.

\begin{prop}[\cite{ZH}, \cite{MaZ.3, HZ}]\label{altGFbounds}
Under the assumptions of Proposition \ref{greenbounds}, $\gtild$ has
the following bounds.
\begin{equation}\label{Gbounds}
\begin{aligned}
|\partial_{x,y}^\alpha &\tilde G(x,t;y)|\le \\
& C(t^{-|\alpha|/2} + |\alpha_x| e^{-\eta|x|}) \Big( \sum_{k=1}^n
t^{-1/2}e^{-(x-y-a_k^{-} t)^2/Mt} e^{-\eta x^+} \\
&+ \sum_{a_k^{-} > 0, \, a_j^{-} < 0} \chi_{\{ |a_k^{-} t|\ge |y|
\}} t^{-1/2} e^{-(x-a_j^{-}(t-|y/a_k^{-}|))^2/Mt}
e^{-\eta x^+}, \\
&+ \sum_{a_k^{-} > 0, \, a_j^{+}> 0} \chi_{\{ |a_k^{-} t|\ge |y| \}}
t^{-1/2} e^{-(x-a_j^{+} (t-|y/a_k^{-}|))^2/Mt}
e^{-\eta x^-}\Big), \\
\end{aligned}
\end{equation}
for $y\le 0$, and symmetrically for $y\ge 0$, for some $\eta$, $C$,
$M>0$, where $a_j^\pm$ are as in Proposition \ref{greenbounds},
$\beta_k^\pm>0$, $x^\pm$ denotes the positive/negative part of $x$,
and indicator function $\chi_{\{ |a_k^{-}t|\ge |y| \}}$ is $1$ for
$|a_k^{-}t|\ge |y|$ and $0$ otherwise. Moreover, all estimates are
uniform in the supressed parameter $\delta_*$.
\end{prop}
\begin{rem}\label{GFremarks}
\textup{We will refer to the three differently scaled diffusion
kernels in (\ref{Gbounds}) respectively as the {\it convection}
kernel, the {\it reflection} kernel, and the {\it transmission}
kernel. We recall the notation
$$
\textrm{errfn} (z) := \frac{1}{2\pi} \int_{-\infty}^z e^{-\xi^2}
d\xi.
$$
}
\end{rem}
\begin{rem}
\textup{ Green function bounds for the strictly parabolic case (see
Remark \ref{strictpc}) is very similar. The main difference  is
that, in the strictly parabolic case,  the hyperbolic part $H$ is
absent in the decomposition of $G$ as in (\ref{ourdecomp}). }
\end{rem}
\begin{rem}\label{eboundsrmk}
\textup{ {}From \eqref{e}, we obtain by straightforward calculation
(see \cite{MaZ.3}) the bounds
\begin{equation}\label{ebounds}
\begin{aligned}
|e_j(y,t)|&\le C\sum_{a_k^->0}
  \left(\textrm{errfn }\left(\frac{y+a_k^{-}t}{\sqrt{4\beta_k^{-}t}}\right)
  -\textrm{errfn }\left(\frac{y-a_k^{-}t}{\sqrt{4\beta_k^{-}t}}\right)\right),\\
|\partial_t  e_j(y,t)|&\le C t^{-1/2} \sum_{a_k^->0} e^{-|y+a_k^-t|^2/Mt},\\
|\partial_y  e_j(y,t)|&\le C t^{-1/2} \sum_{a_k^->0} e^{-|y+a_k^-t|^2/Mt}\\
|\partial_{yt}  e_j(y,t)|&\le C
t^{-1} \sum_{a_k^->0} e^{-|y+a_k^-t|^2/Mt}\\
\end{aligned}
\end{equation}
for $y\le 0$, and symmetrically for $y\ge 0$. }
\end{rem}

\section{Nonlinear analysis}

Let $\utild$ solve (\ref{viscous}), and, using ($\dcal$ii), assume
that
$$\int_{-\infty}^{+\infty} \utild (x, 0) - \ubar(x)= \sum_{a_j^- <0}m_j r_j^- +
\sum_{a_j^+ >0}m_j r_j^+ +\sum_{i=1}^{\ell}\int c_i\frac{\partial
\ubar^\delta}{\partial \delta_i}_{|\delta=0}
$$
with $m_i$'s and $c_i$'s small enough. Using the Implicit Function
Theorem, we can find $\delta_*$ such that
$$\int_{-\infty}^{+\infty} \utild (x, 0) - \ubar^{\delta_*}(x)= \sum_{a_j^- <0}m'_j r_j^- +
\sum_{a_j^+ >0}m'_j r_j^+
$$
where each $m'_i$ is just ``slightly" different from $m_i$. Notice
that this way we have no ``mass" in any $\int \frac{\partial
\ubar^\delta}{\partial \delta_i}$ direction anymore. With a slight
abuse of notation we drop the $'$ sign from each $m'_i$ and denote
it simply by $m_i.$
\begin{rem}\textup{
In the case of Lax--type shock waves,
$\ubar^\delta(x)=\ubar(x+\delta),$ hence $\ddp=u'(x),$ and
$\delta_*$ can be explicitly computed: $\delta_*=c_1.$ }
\end{rem}
Let $u(x,t)=\utild(x,t)-\ubar^{\delta_*}(x)$  and use Taylor's
expansion around $\ubar^{\delta_*}(x)$ to find
\begin{equation}
u_t + (A(x)u)_x - (B(x)u_x)_x = -(\Gamma(x)(u,u))_x + Q(u, u_x)_x,
\label{utaylor}
\end{equation}
where $\Gamma(x)(u,u) = d^2f(\ubar^{\delta_*})(u,u)-
d^2B(\ubar^{\delta_*})(u,u)\ubar^{\delta_*}_x$ and $$Q(u, u_x) =
\mathbf {O} (|u||u_x|+|u|^3).$$  Denote $\Gamma^\pm = \Gamma(\pm
\infty).$  Define constant coefficients $b^{\pm}_{ij}$ and
$\Gamma_{ijk}^{\pm}$ to satisfy
\begin{equation}
\Gamma^\pm (r^\pm_j, r^\pm_k) = \sum_{i=1}^n \Gamma_{ijk}^{\pm}
r^\pm_i,\quad B^\pm r^\pm_j= \sum_{i=1}^n b^{\pm}_{ij} r^\pm_i.
\label{coefshock}
\end{equation}
Then, of course, $\beta^\pm_i = b^{\pm}_{ii}$ and $\gamma^\pm_i :=
\Gamma_{iii}^{\pm}.$

Now define $\fe_i^-$ by (\ref{burgers}) and (\ref{burgersdata}), and
likewise for $\fe_i^+$. Finally $\fe$ is defined by (\ref{phi}).
Then set $v:= u -\fe-\ddp\delta(t)$, with $\delta(t)$ to be defined
later, and assuming $\delta(0)=0$. Notice that, by our choice of
$\delta_*$ and diffusion waves $\fe_i^\pm$'s, we have zero initial
mass of $v$, i.e.,
\begin{equation} \label{initial0mass}
\int_{-\infty}^{+\infty} v(x,0)dx = 0
\end{equation}

Replacing $u$ with $v + \fe+\ddp\delta(t)$ in (\ref{utaylor}) (
$\frac{\partial \ubar^\delta}{\partial \delta_i} $  computed at
$\delta=\delta_*$), and using the fact that
$\frac{\partial\ubar^\delta}{\partial \delta_i}$ satisfies the
linear time independent equation $Lv=0$, we will have
\begin{equation}
v_t - Lv = \Phi(x,t)+ \fcal(\fe, v,  \ddp\delta(t))_x + \ddp\dot
\delta (t), \label{khati}
\end{equation}
where
\begin{equation}\label{fcalt}
\begin{split}
\fcal(\fe, v, \ddp\delta &) = \mathbf {O}  ( |v|^2 + |\fe||v|
+ |v||\ddp\delta|+ |\fe|| \ddp\delta|+|\ddp\delta|^2 \\
&+ |(\fe + v+\ddp\delta)(\fe + v + \ddp\delta)_x| + |\fe
+v+\ddp\delta|^3 ).
\end{split}
\end{equation}
Furthermore
\begin{equation}
\begin{aligned}
\fcal(v, \fe, \ddp\delta(t))_x& =
\mathbf{O}\big(\fcal(v,\fe,\ddp\delta)\\
&+
|(v+\fe+\ddp\delta)_x||(v^{II}+\fe+\ddp\delta)_x|\\
&+|v+\fe+\ddp\delta||(v^{II}+\fe+\ddp\delta)_{xx}|\big),
\label{fcalx}
\end{aligned}
\end{equation}
and $\Phi(x,t) := -\fe_t - (A(x)\fe)_x  + (B(x)\fe_x)_x
-(\Gamma(x)(\fe , \fe))_x$. For $\Phi$  we write
\begin{equation}
\begin{split}
\Phi(x,t) = &-(\fe_t +A \fe_x - B \fe_{xx} + \Gamma(\fe, \fe)_x) \\
= &-\sum_{a_i^- < 0} \fe_t^i r_i^- + (A(x)\fe^i r_i^-)_x -
(B(x)\fe_x^i r_i^-)_x + (\Gamma(x)(\fe^ir_i^-,\fe^ir_i^-))_x\\
&- \sum_{a_i^+ > 0} \fe_t^i r_i^+ + (A(x)\fe^i r_i^+)_x -
(B(x)\fe_x^i r_i^+)_x + (\Gamma(x)(\fe^ir_i^+,\fe^ir_i^+))_x\\
& - \sum_{i\neq j}(\fe_i \fe_j \Gamma(x)(r_i^\pm, r_j^\pm))_x.
\end{split} \label{ghati}
\end{equation}
Let us write a typical term of the first summation ($a_i^- < 0$)in
the following form:
\begin{equation}
\begin{split}
\fe_t^i & r_i^- + (A(x)\fe^i r_i^-)_x -
(B(x)\fe_x^i r_i^-)_x + (\Gamma(x)(\fe^ir_i^-,\fe^ir_i^-))_x\\
&= \left[(A(x)-A^-)\fe^i r_i^- -
(B(x)-B^-)\fe_x^i r_i^- + (\Gamma(x)-\Gamma^-)(\fe^ir_i^-,\fe^ir_i^-)\right]_x\\
&+ \fe_t^i r_i^- + (\fe^i_xA^- r_i^-) - (\fe_{xx}^i B^- r_i^-) +
((\fe^{i})^2_x\Gamma^-(r_i^-,r_i^-)).
\end{split} \label{ghatitar}
\end{equation}
Now we use the definition of $\fe^i$ in (\ref{burgers}) and the
definition of coefficients $b_{ij}$ and $\Gamma_{ijk}$ in
(\ref{coefshock}) to write the last part of (\ref{ghatitar}) in the
following form:
\begin{equation}
\begin{split}
\fe_t^i r_i^- + (\fe^i_xA^- r_i^-) - &(\fe_{xx}^i B^- r_i^-) +
((\fe^i)^2_x\Gamma^-(r_i^-,r_i^-))\\
&= -\fe^i_{xx}\sum_{j\ne i} b^-_{ij}r^-_j - (\fe^i)^2_{x}\sum_{j\ne
i} \Gamma^-_{jii}r^-_j.
\end{split} \label{ajabaa}
\end{equation}
Similar statements hold for $a_i^+ > 0$ with minus signs replaced
with plus signs.

\medskip
Now we employ Duhamel's principle to get from (\ref{khati})
\begin{equation}
\begin{aligned} v(x,t)
&=\int^{+\infty}_{-\infty} G(x,t;y)v_0(y)dy\\
&+\int^t_0
\int^{+\infty}_{-\infty}G(x,t-s;y)\fcal(\fe, v, \ddp\delta)_y(y,s) dy \, ds,\\
&+ \int^t_0\int^{+\infty}_{-\infty}G(x,t-s;y)\Phi(y,s) dy
\, ds \\
&+ \ddp \delta(t)\\
 \label{shswgreen}
\end{aligned}
\end{equation}
To obtain the above, we used
$$\int_{-\infty}^{+\infty}
G(x,t;y)\frac{\partial\ubar^\delta}{\partial
\delta_i}(y)dy=e^{Lt}\frac{\partial\ubar^\delta}{\partial
\delta_i}=\frac{\partial\ubar^\delta}{\partial \delta_i}$$ and
$\delta(0)=0$.

Set
\begin{equation}
\begin{aligned}
 \delta_i (t)
&=-\int^\infty_{-\infty}e_{i}(y,t) v_0(y)dy\\
&-\int^t_0\int^{+\infty}_{-\infty} e_{i}(y,t-s)\fcal(\fe,
v,\ddp\delta)_y(y,s)
 dy ds\\
&-\int^t_0\int^{+\infty}_{-\infty}e_i(x,t-s;y)\Phi(y,s)dy \, ds.
\end{aligned} \label{delta}
\end{equation}
Using (\ref{shswgreen}),  (\ref{delta}) and $G=H+E+\gtild$ we
obtain:
\begin{equation}
\begin{aligned} v(x,t)
&=\int^{+\infty}_{-\infty} (H+ \gtild)(x,t;y)v_0(y)dy\\
&+\int^t_0
\int^{+\infty}_{-\infty}(H+\gtild)(x,t-s;y)\fcal(\fe, v,\ddp\delta)_y(y,s) dy \, ds\\
&+ \int^t_0\int^{+\infty}_{-\infty}(H+\gtild)(x,t-s;y)\Phi(y,s) dy
\, ds.  \label{gtild}
\end{aligned}
\end{equation}

In addition to $v(x, t)$ and $\delta (t)$, we will keep track in our
argument of $v_x (x, t)$ and $\dot \delta (t)$, the latter of which
satisfies
\begin{equation}
\begin{aligned}
 \dot \delta_i (t)
&=-\int^\infty_{-\infty}\partial_t e_{i}(y,t) v_0(y)dy\\
&-\int^t_0\int^{+\infty}_{-\infty} \partial_t e_{i}(y,t-s)\fcal(\fe,
v,\ddp\delta)_y(y,s)
 dy ds\\
&-\int^t_0\int^{+\infty}_{-\infty} \partial_t e_i(y,t-s)\Phi(y,s)dy \, ds,
\end{aligned} \label{dotdelta}
\end{equation}
where we have taken advantage of the observation, apparent
from (\ref{ebounds}), that $e(y,0)=0$.

The following lemmas are the main ingredients  for the proof of
Theorem \ref{thmmain}. Their proofs, however, are postponed to
Section \ref{S:integralEst}.

\begin{lem}[Estimates for linear part] \label{linear}
Suppose that for some $E_0>0$, we have that $v_0 (y)$ satisfies the
conditions
\begin{equation*}
\begin{aligned}
|v_0 (y)| &\le E_0 (1 + |y|)^{-3/2} \\
\int_{-\infty}^{+\infty} &v_0 (y) dy = 0.
\end{aligned}
\end{equation*}
Then there holds
\begin{align}
\Big|\int_{-\infty}^{+\infty} \tilde{G}(x,t;y) v_0 (y) dy \Big| &
\le
C E_0 \psi_2 (x,t).\label{gv0} \\
\Big|\int_{-\infty}^{+\infty} e_i (y,t) v_0 (y) dy \Big| & \le
C E_0 (1+t)^{-1/2},\label{ev0} \\
\Big|\int_{-\infty}^{+\infty} \partial_t e_i (y,t) v_0 (y) dy \Big|
& \le C E_0 (1+t)^{-3/2},\label{etv0}
\end{align}
where $C$ does not depend on $E_0$.
\end{lem}
\begin{lem}\label{Hlinear}
If $|v_0(x)|\le E_0 (1+|x|)^{-\frac32},$  then
\begin{equation}
\int_{-\infty}^{+\infty} H(x, t; y) v_0(y) dy\le C E_0 e^{-\theta
t}(1+|x|)^{-\frac32} \label{Hv0}
\end{equation}
If $|\partial_x v_0(x)|\le E_0 (1+|x|)^{-\frac12},$ then
\begin{equation}
\int_{-\infty}^{+\infty} H_x(x, t; y) v_0(y) dy\le C E_0 e^{-\theta
t}(1+|x|)^{-\frac12}\label{Hxv0}
\end{equation}
for some $\theta > 0$. The right hand sides of (\ref{Hv0}) and
(\ref{Hxv0}) are obviously of  order $\alpha(x,t).$
\end{lem}
\medskip
\begin{lem}[Estimates for nonlinear part]\label{nonlinearestimates}
For $G(x,t;y)$ as in Propositions \ref{greenbounds} and
\ref{altGFbounds}, we have
\begin{align}
\int_0^t \int_{-\infty}^{+\infty} | \tilde{G}_y (x,t-s;y)| \Psi
(y,s) dyds &\le
C (\bar\psi_1 + \psi_2 + \alpha) (x,t), \label{gpsi}\\
\Big| \int_0^t \int_{-\infty}^{+\infty}  \tilde{G} (x,t-s;y) \Phi
(y,s) dyds \Big| &\le
CE_0 (\bar\psi_1 + \psi_2) (x,t), \label{gphi}\\
\int_0^t \int_{-\infty}^{+\infty} | \partial_y e_i (y,t-s)| \Psi
(y,s) dyds &\le
C (1+t)^{-3/4},\label{epsi} \\
\Big| \int_0^t \int_{-\infty}^{+\infty}  e_i (y,t-s) \Phi (y,s) dyds
\Big| &\le
CE_0 (1+t)^{-1/2}, \label{ephi}\\
\Big| \int_0^t \int_{-\infty}^{+\infty}  \partial_t e_i (y,t-s) \Phi
(y,s) dyds \Big| &\le
C E_0(1+t)^{-1}, \label{etphi}\\
\int_0^t \int_{-\infty}^{+\infty} | \partial_{yt} e_i (y,t-s)| \Psi
(y,s) dyds &\le C (1+t)^{-1},\label{etpsi}
\end{align}
for $\Phi(y, s)$ as in (\ref{ghati}) and
\begin{equation}
\begin{aligned}
\Psi(y,s) &= (1+s)^{-1/4} s^{-1/2} (\bar\psi_1 + \psi_2 + \alpha + \varphi) (y,s) \\
& + (1+s)^{-1/2} s^{-1/2} e^{-\eta |y|}.
\end{aligned} \label{Psi}
\end{equation}
\end{lem}
\medskip
\begin{lem}\label{Hnonlinear} If  \,$|\Upsilon(y,s)|\le s^{-1/2}(\bar\psi_1 +
\psi_2 + \alpha + \varphi )(y,s)+ s^{-1/2}e^{-\eta |y|},$ then
\begin{equation}
 \Big|\int_0^t \int_{-\infty}^{+\infty}  H (x,t-s;y)  \Upsilon(y,s) dyds\Big| \le C ( \bar\psi_1 + \psi_2 +
\alpha) (x,t). \label{Hf}
\end{equation}
If $|\partial_y \Upsilon(y,s)|\le  s^{-1/2}(\bar\psi_1 + \psi_2 +
\alpha + \varphi + e^{-\eta |y|})(y,s),$ then
\begin{equation}
 \Big|\int_0^t \int_{-\infty}^{+\infty}  H_x (x,t-s;y)  \Upsilon(y,s) dyds\Big| \le C ( \bar\psi_1 + \psi_2 +
\alpha) (x,t). \label{Hxf}
\end{equation}
Also,
\begin{equation}
 \Big|\int_0^t \int_{-\infty}^{+\infty}  H (x,t-s;y)  \Phi(y,s) dyds\Big| \le C E_0 ( \bar\psi_1 + \psi_2 +
\alpha) (x,t)  (x,t), \label{Hphi}
\end{equation}
and similarly with $H$ replaced by $H_x.$
\end{lem}

\medskip
\noindent {\bf Proof of Theorem \ref{thmmain}.} We use estimates
(\ref{vlp}) and (\ref{vhs})  in order to obtain a supremum  norm on
$v$ and its first,  second and third derivatives in $x$. We prove
(\ref{1stestimate}) with $\psi_1$ replaced by $\bar\psi_1$,i.e. we
prove
\begin{equation}
|v(x,t)|\le
 C E_0
(\bar\psi_1+\psi_2 +\alpha); \label{barpsiest},
\end{equation}
 but then, as we have $|v(x,t)|\le CE_0(1+t)^{-\frac34}$ (this
is (\ref{vlp}) for $p=\infty$), and as
$$\min\{(1+t)^{-\frac34}, \, \bar\psi_1(x,t)\}\sim \psi_1(x,t),$$
the proof of (\ref{1stestimate}) would be immediate.

\smallskip
Let
\begin{equation}
\begin{aligned}
    \zeta(t) := \sup_{y, 0\le s \le t}\frac{|v(y,
    s)|}{\big(\bar\psi_1+\psi_2+\alpha \big)(y,s)}
    + \sup_{y, 0\le s \le t}\frac{|v_x(y,
    s)|}{t^{-\frac12}(1+t)^{\frac12}\big(\bar\psi_1+\psi_2+\alpha\big)(y,s)}\\
    + \sup_{0\le s \le t} |\delta (s)|(1+s)^{\frac 12 }
     + \sup_{0\le s
    \le t} |\dot\delta (s)|(1+s).
\end{aligned} \label{zetadef}
\end{equation}
Our aim is to show that
\begin{equation}
       \zeta (t) \leq C(E_0 + E_0\zeta(t)), \label{zetaineq}
\end{equation}
from which we conclude $\zeta(t)\le \frac{CE_0}{1-CE_0}\le \frac12$,
 if $E_0\le \frac1{2 C}$.
Equivalent to (\ref{zetaineq}) is
\begin{equation}
      |v(x, t)| \leq C(E_0
       +E_0\zeta(t)) (\bar\psi_1+\psi_2+\alpha)(y,s),
\end{equation}
\begin{equation}
    |\delta(t)| \leq C(E_0
       +E_0\zeta(t)) (1+t)^{-\frac 12 },
\label{delbd}
\end{equation}
and similar statements for  $v_x$ and $\dot\delta$.
\medskip\\
\emph{Estimates of $v(x,t)$.} Looking at (\ref{gtild}), there are
three parts that should be estimated. The first part, the linear
part, is carried out using (\ref{Hv0}) and (\ref{gv0}). For the
second part, we notice that, by  (\ref{fcalt}), (\ref{zetadef}),
(\ref{vlp}), (\ref{deltabound}),  (\ref{dotdeltabd}), Lemma
\ref{expdecay}, and the definition and bounds of $\fe,$
$$|\fcal(y,s)|\le C(E_0+\zeta(t))\Psi(y,s)$$
with $\fcal(y,s)$ as in (\ref{fcalt}), and $\Psi$ as in (\ref{Psi}).
Hence, by integration by part,
\begin{equation*}
\begin{aligned}
\Big|\int^t_0 &\int^{+\infty}_{-\infty}\gtild(x,t-s;y)\fcal(\fe,
v,\ddp\delta)_y(y,s) dy \, ds \Big|\\
&= \Big| \int^t_0
\int^{+\infty}_{-\infty}\gtild_y(x,t-s;y)\fcal(\fe,
v,\ddp\delta)(y,s) dy \, ds\Big|\\
&\leq C(E_0
       +E_0\zeta(t)) h(x,t),
\end{aligned}
\end{equation*}
by (\ref{gpsi}).

 As we do not have good estimates for $H_y$ we cannot do the same
with part containing $H$, so, instead, we notice that
$$|\fcal_y(y,s)|\le C(E_0+\zeta(t))\Upsilon(y,s)$$
$\fcal_y(y,s)$ as in (\ref{fcalx}), and $\Upsilon(y,s)$ as in Lemma
\ref{Hnonlinear}. (The subtle fact is that, as we have $v_{xx}$ in
$\fcal_y(y,s)$, we need to use (\ref{vhs}), hence having only
$(1+t)^{-\frac12}$, instead of $(1+t)^{-\frac34}$ which appears in
$\Psi$.) Therefore we can use (\ref{Hf}) to obtain the desired
estimate. The third integral of (\ref{gtild}) can be estimated
similarly using (\ref{gphi}) and (\ref{Hphi}).
\medskip\\
\emph{Estimates of $\delta(t)$ and $\dot\delta(t)$.} Similarly and
using (\ref{delta}) and (\ref{epsi}--(\ref{etpsi}).
\medskip\\
\emph{Estimates of $v_x(x,t)$.} By (\ref{gtild}),
\begin{equation}
\begin{aligned} v_x(x,t)
&=\int^{+\infty}_{-\infty} (H_x+ \gtild_x)(x,t;y)v_0(y)dy\\
&+\int^t_0
\int^{+\infty}_{-\infty}(H_x+\gtild_x)(x,t-s;y)\fcal(\fe, v,\ddp\delta)_y(y,s) dy \, ds,\\
&+ \int^t_0\int^{+\infty}_{-\infty}(H_x+\gtild_x)(x,t-s;y)\Phi(y,s)
dy \, ds.  \label{vxform}
\end{aligned}
\end{equation}
The parts involving $H_x$ are treated using (\ref{Hxv0}),
(\ref{Hxf}) and the similar statement for $\Phi$ in Lemma
\ref{Hnonlinear}; this because, due to (\ref{vhs}) $\fcal_{yy}(y,s)$
has the bounds similar to $\fcal_y(y,s).$ For
$\int^{+\infty}_{-\infty} \gtild_x (x,t;y)v_0(y)dy$, we notice that
$\gtild_x$ is at least as good as $\gtild$  away from $0$.
Therefore, for $t\ge 1$, a proof identical to that of (\ref{gv0})
can also prove
$$\int^{+\infty}_{-\infty} \gtild_x
(x,t;y)v_0(y)dy\le CE_0 \psi_2(x,t).$$ On the other hand, for $t\le
1$ we use the fact that $\gtild_x(x,t;y \sim t^{-\frac12}
\gtild(x,mt;y)$ for $t$ near zero and some positive $m$. This,
together with (\ref{vhs}) provides us with the necessary bounds for
the linear part of the calculations.

Finally,
\begin{equation*}
\begin{aligned}
\int^t_0
\int^{+\infty}_{-\infty}&\gtild_x(x,t-s;y)\fcal(\fe, v,\ddp\delta)_y(y,s) dy \, ds\\
&= \int_0^{t-1}\int^{+\infty}_{-\infty}\gtild_{xy}(x,t-s;y)\fcal(\fe, v,\ddp\delta)(y,s) dy \, ds\\
&+\int_{t-1}^t\int^{+\infty}_{-\infty}\gtild_x(x,t-s;y)\fcal(\fe,
v,\ddp\delta)_y(y,s) dy \, ds.
\end{aligned}
\end{equation*}
The first integral can be computed exactly the same way (\ref{gpsi})
is proved, for away from $t-s=0$, $\gtild_{xy}(x,t-s;y)$ is as well
as $\gtild_{xy}(x,t-s;y)$. For the second integral, a close
investigation of the proof of (\ref{gpsi}) shows us that, if we
limit the time integration to $t-1\le s \le t$, then the same proof
works to prove that
$$\int_{t-1}^t\int_{-\infty}^{+\infty} | \tilde{G}_x (x,t-s;y)| |\Upsilon (y,s)| dyds
\le C (\bar\psi_1 + \psi_2 + \alpha) (x,t),$$ with $\Upsilon$ as in
Lemma \ref{Hnonlinear}. This fact plus the bound for $\fcal_y$
 provides us with the estimates needed for the second integral.

 This concludes the proof of Theorem \ref{thmmain}. \hfill $\square$

\section{Integral estimates}\label{S:integralEst}
\noindent {\bf Proof of Lemma \ref{linear}.} We begin by
establishing that for $\tilde{G}$ satisfying
\begin{equation*}
\|\tilde{G}\|_{L^1_x(\mathbb{R})} \le C_1, \quad
\|\tilde{G}_y\|_{L^\infty_x(\mathbb{R})} \le C_2 t^{-1}, \quad
\|\tilde{G}_y\|_{L^1_x(\mathbb{R})} \le C_3 t^{-1/2},
\end{equation*}
for some positive constants $C_1$, $C_2$, and $C_3$, we have
\begin{equation} \label{mintime}
\Big|\int_{-\infty}^{+\infty} \tilde{G}(x,t;y) v_0 (y) dy \Big| \le
C E_0 (1+t)^{-3/4}.
\end{equation}
First, in the event that $t$ is bounded, we have
\begin{equation*}
\Big|\int_{-\infty}^{+\infty} \tilde{G}(x,t;y) v_0 (y) dy \Big| \le
\|\tilde{G}\|_{L^1_x(\mathbb{R})} \|v_0\|_{L^\infty} \le C_1 E_0.
\end{equation*}
In the alternative case, for which we take $t$ bounded away from 0,
we integrate by parts to obtain
\begin{equation*}
\int_{-\infty}^{+\infty} \tilde{G}(x,t;y) v_0 (y) dy =
\int_{-\infty}^{+\infty} \tilde{G}_y(x,t;y) V_0 (y) dy,
\end{equation*}
where
\begin{equation*}
V_0(x):= \int_{-\infty}^x v_0(x) dx, \qquad |V_0 (y)| \le C_4 E_0
(1+|y|)^{-1/2}.
\end{equation*}
We have, then,
\begin{equation*}
\begin{aligned}
\Big|& \int_{-\infty}^{+\infty} \tilde{G}_y (x,t;y) V_0 (y) dy \Big| \\
&\le \Big|\int_{\{|y|\le\sqrt{t} \}} \tilde{G}_y (x,t;y) V_0 (y) dy
\Big| +
\Big|\int_{\{|y|\ge\sqrt{t} \}} \tilde{G}_y (x,t;y) V_0 (y) dy \Big| \\
& \le C_2 t^{-1} \int_{\{|y|\le\sqrt{t} \}} |V_0 (y)| dy +
C_4 E_0 (1+\sqrt{t})^{-1/2} \int_{\{|y|\ge\sqrt{t} \}} |\tilde{G}_y (x,t;y)| dy  \\
&\le C_5 E_0 t^{-3/4},
\end{aligned}
\end{equation*}
establishing (\ref{mintime}).

It now follows from (\ref{mintime}) that if there exists some
$a_j^\pm \gtrless 0$ so that $|x - a_j^\pm t| \le \sqrt{t}$, then
$t^{-3/4} \le C|x-a_j^\pm t|^{-3/2}$, and the first estimate is
apparent. In the event that $|x - a_j^{\pm} t| \ge \sqrt{t}$ for all
$a_j^\pm$, we consider the case in which there exists some $a_j^\pm
\gtrless 0$ so that $|x - a_j^\pm t| \le \epsilon t$, where
$\epsilon > 0$ is sufficiently small so that for all $j \ne k$, $|x
- a_k^\pm t| \ge \eta t$, for some $\eta > 0$.  In this case, we
compute
\begin{equation} \label{linearparts}
\begin{aligned}
\int_{-\infty}^{+\infty}& \tilde{G} (x,t;y) v_0 (y) dy \\
&= \int_{\{|y|\le \frac{|x - a_j^\pm t|}{N} \}} \tilde{G} (x,t;y)
v_0 (y) dy + \int_{\{|y|\ge \frac{|x - a_j^\pm t|}{N} \}} \tilde{G}
(x,t;y) v_0 (y) dy,
\end{aligned}
\end{equation}
where $N$ will be chosen sufficiently large in the analysis.  For
the second of these last two integrals, integration of $\tilde{G}$
immediately gives an estimate by
\begin{equation*}
C E_0 (1 + |x - a_j^\pm t|)^{-3/2},
\end{equation*}
which is sufficient for our first estimate since we are in the case
$|x - a_j^\pm t| \ge \sqrt{t}$.  Alternatively, for the first
integral in (\ref{linearparts}), we integrate by parts to obtain
\begin{equation} \label{linearcases}
\begin{aligned}
\Big| &\int_{\{|y|\le \frac{|x - a_j^\pm t|}{\bar{L}} \}} \tilde{G} (x,t;y) v_0 (y) dy \Big| \\
&\le \Big|\tilde{G} (x,t,\frac{x - a_j^\pm t}{N}) V_0 (\frac{x -
a_j^\pm t}{N})  \Big| +
\Big|\tilde{G} (x,t,-\frac{x - a_j^\pm t}{N}) V_0 (-\frac{x - a_j^\pm t}{N})  \Big| \\
&+ \int_{\{|y|\le \frac{|x - a_j^\pm t|}{N} \}} |\tilde{G}_y
(x,t;y)| |V_0 (y)| dy.
\end{aligned}
\end{equation}
It remains to estimate each of these last three terms for each
summand in the estimates on $|\tilde{G}|$ and $|\tilde{G}_y|$ from
Lemma \ref{altGFbounds}.  As each case is similar, we proceed only
with estimates on the convection terms.

In the case $y < 0$ and for the convection Green's function
estimate, we have
\begin{equation*}
\begin{aligned}
&\Big|\tilde{G}(x,t;-\frac{|x-a_j^- t|}{N}) V_0 (-\frac{|x - a_j^- t|}{N})\Big| \\
&\le CE_0 \sum_{k = 1}^n t^{-\frac{1}{2}} e^{-\frac{(x -
\frac{|x-a_j^- t|}{N}-a_k^- t)^2}{Mt}} e^{-\eta x^+} (1 + |x - a_j^-
t|)^{-\frac{1}{2}}.
\end{aligned}
\end{equation*}
In the event that $j=k$, we have
\begin{equation*}
|x - a_k^- t| - (\frac{x - a_k^- t}{N}) \ge (x - a_k^- t) (1 -
\frac{1}{N}),
\end{equation*}
which for $N$ sufficiently large leads immediately to an estimate by
\begin{equation*}
C E_0 t^{-1/2} e^{-\frac{(x - a_k^- t)^2}{Lt}} (1 + |x - a_k^-
t|)^{-1/2},
\end{equation*}
for some $L$ sufficiently large.  According to the boundedness of
\begin{equation} \label{kerneltotail}
\frac{|x - a_k^- t|}{t^{1/2}} e^{-\frac{(x - a_k^- t)^2}{Lt}},
\end{equation}
we obtain the claimed estimate.  On the other hand, if $j \ne k$,
then $|x - a_j^- t| \le \epsilon t$ and $|x - a_k^- t| \ge \eta t$,
so that for $N$ sufficiently large, we have exponential decay in
$t$.

We next consider the integral in (\ref{linearcases}), for which in
the case of the convection estimate and for $y \le 0$, we have
\begin{equation*}
\int_{\{y \le -\frac{|x - a_j^- t|}{N}\}} \sum_{k=1}^n t^{-1}
e^{-\frac{(x - \frac{|x-a_j^- t|}{N}-a_k^- t)^2}{Mt}} e^{-\eta x^+}
|V_0 (y)| dy.
\end{equation*}
Proceeding as in the boundary case, we see that in the case $j \ne
k$, the kernel decays at exponential rate in $t$, while for $j = k$,
upon integration of $V_0 (y)$, we have an estimate,
\begin{equation*}
\begin{aligned}
C &E_0 t^{-1} e^{-\frac{(x - a_j^- t)^2}{Lt}} (1 + |x - a_j^- t|)^{1/2} \\
&\le \bar{C} E_0 |x - a_j^- t|^{-3/2}.
\end{aligned}
\end{equation*}

For the second estimate in Lemma \ref{linear}, we observe from
Remark \ref{eboundsrmk} the estimate
\begin{equation*}
|\partial_y e_i (y,t)| \le C t^{-1/2} \sum_{a_k^- > 0}
e^{-\frac{(y+a_k^- t)^2}{Mt}}.
\end{equation*}
Integrating by parts, then, we have
\begin{equation*}
\begin{aligned}
\Big|&\int_{-\infty}^{+\infty} e_i (y,t) v_0 (y) dy \Big|
= \Big|\int_{-\infty}^{+\infty} \partial_y e_i (y,t) V_0 (y) dy \Big| \\
&\le C t^{-1/2} \sum_{a_k^- > 0} \int_{-\infty}^{+\infty}
e^{-\frac{(y+a_k^- t)^2}{Mt}} E_0 (1+|y|)^{-1/2} dy.
\end{aligned}
\end{equation*}
In this last integrand, we observe the inequality
\begin{equation*}
e^{-\frac{(y+a_k^- t)^2}{Mt}} E_0 (1+|y|)^{-1/2} \le C
e^{-\frac{(y+a_k^- t)^2}{\bar{M}t}} E_0 (1+t)^{-1/2},
\end{equation*}
from which the claimed estimate is immediate.

The final estimate in Lemma \ref{linear} is proven similarly as the
second estimate in Lemma 2 of \cite{HZ}. \hfill $\square$

\medskip
\noindent {\bf Proof of Lemma \ref{Hlinear}.} Looking at
(\ref{multH}), we notice that in order to estimate
$\int_{-\infty}^{+\infty}H(x,t;y)v_0(y) dy$ it suffices to estimate
$$\int_{-\infty}^{+\infty}\mathcal{R}_j^*(x) \mathcal{O}(e^{-\eta_0 t}) \delta_{x-\bar a_j^*
t}(-y) \mathcal{L}_j^{*t}(y)v_0(y) dy$$
$$\le CE_0 e^{-\eta_0 t}v_0(\bar a_j^*t-x)$$
$$\le CE_0 e^{-\eta_0 t} (1+|\bar a_j^*t-x|)^{-\frac32} $$
$$\le CE_0 e^{-\eta_0 t} (1+|x|)^{-\frac32}(1+|\bar a_j^*t|)^{\frac32}$$
$$\le CE_0 e^{-\frac{\eta_0 t}2} (1+|x|)^{-\frac32}$$
Here we used the crude inequality
\begin{equation}\label{abineq}
\frac {1}{1+|a+b|}\le \frac{1+|b|}{1+|a|}
\end{equation}
and the fact that $\bar a_j^*$ $\mathcal{R}_j^*$ and
$\mathcal{L}_j^{*t}$ are bounded. This gives us (\ref{Hv0}).
Estimate (\ref{Hxv0}) is obtained similarly.
 \hfill $\square$

\medskip
\noindent {\bf Proof of Lemma \ref{nonlinearestimates}.} For Lemma
\ref{nonlinearestimates}, the proof of each estimate requires the
analysis of several cases.  We proceed by carrying out detailed
calculations in the most delicate cases and sufficing to indicate
the appropriate arguments in the others.  In particular, we will
always consider the case $x, y \le 0$.  The case $y \le 0 \le x$ is
similar (though certainly not identical) to the reflection case for
$x, y \le 0$, and the estimates for $y \ge 0$ are entirely
symmetric.  The analysis of each type of kernel---convection,
transmission, and scattering---is similar, and we carry out details
only in the case of convection.  (This terminology is reviewed in
Remark \ref{GFremarks}).

{\it Nonlinearity} $(1+s)^{-1/4} s^{-1/2} \bar\psi_1$.  We begin by
estimating integrals,
\begin{equation} \label{case1integral}
\int_0^t \int_{-|a_1^-|s}^0 (t-s)^{-1} e^{-\frac{(x - y - a_j^-
(t-s))^2}{M(t-s)}} (1 + |y - a_k^- s|)^{-3/4} s^{-1/2} (1+s)^{-3/4}
dy ds.
\end{equation}
In the event that $|x| \ge |a_1^-| t$, we write
\begin{equation} \label{firstdecomp}
x - y - a_j^- (t-s) = (x - a_1^- t) - (y - a_1^- s) + (a_1^- -
a_j^-) (t-s),
\end{equation}
and observe that in the current setting ($x \le 0$, $y \in [-|a_1^-|
s, 0]$, $a_1^- \le a_k^-$), there is no cancellation between these
three summands.  Integrating $(1+|y-a_k^- s|)^{-3/4}$ for $s \in [0,
t/2]$ and integrating the kernel for $s \in [t/2,t]$, we obtain an
estimate by
\begin{equation}
\begin{aligned} \label{cancellationargument}
C_1 &t^{-1} e^{-\frac{(x-a_1^- t)^2}{Lt}} \int_0^{t/2} s^{-1/2} (1+s)^{-1/2} ds \\
& +
C_2 t^{-1/2} (1+t)^{-1} e^{-\frac{(x-a_1^- t)^2}{Lt}} \int_{t/2}^t (t-s)^{-1/2} ds \\
&\le C t^{-1} \ln (e+t) e^{-\frac{(x-a_1^- t)^2}{Lt}},
\end{aligned}
\end{equation}
which is sufficient by the argument of (\ref{kerneltotail}).  We
observe that the seeming blowup as $t \to 0$ can be eliminated.
Integrating the kernel in (\ref{case1integral}), we have an estimate
by
\begin{equation} \label{smalltime}
\begin{aligned}
C &\int_0^t (t-s)^{-1/2} s^{-1/2} ds \\
&C_1 t^{-1/2} \int_0^{t/2} s^{-1/2} ds + C_2 t^{-1/2} \int_{t/2}^t
(t-s)^{-1/2} ds \le C.
\end{aligned}
\end{equation}

For $|x| \le |a_1^-| t$, we write
\begin{equation} \label{standarddecomposition}
x - y - a_j^- (t-s) = (x - a_j^- (t-s) - a_k^- s) - (y - a_k^- s),
\end{equation}
from which we observe the inequality
\begin{equation} \label{case1balance}
\begin{aligned}
&e^{-\frac{(x - y - a_j^- (t-s))^2}{M(t-s)}} (1 + |y - a_k^- s|)^{-3/4} \\
&\le C\Big[ e^{-\frac{(x - a_j^- (t-s) - a_k^- s)^2}{\bar{M}(t-s)}}
(1 + |y - a_k^- s|)^{-3/4}
e^{-\epsilon \frac{(x - y - a_j^- (t-s))^2}{M(t-s)}} \\
&+ e^{-\frac{(x - y - a_j^- (t-s))^2}{M(t-s)}} (1 + |y - a_k^- s| +
|x - a_j^- (t-s) -  a_k^- s|)^{-3/4} \Big],
\end{aligned}
\end{equation}
for some $\bar{M}$ sufficiently large and $\epsilon > 0$
sufficiently small. The analysis can now be divided into three
cases: (i) $a_k^- < 0 < a_j^-$, (ii) $a_k^- < a_j^- < 0$, and (iii)
$a_j^- \le a_k^- < 0$.  Each critical argument appears in Case (ii),
and so we consider only it in detail.  For the first estimate in
(\ref{case1balance}), we have
\begin{equation*}
\begin{aligned}
\int_0^t &\int_{-|a_1^-|s}^0 (t-s)^{-1} e^{-\epsilon \frac{(x - y -
a_j^- (t-s))^2}{M(t-s)}}
e^{-\frac{(x - a_j^- (t-s) - a_k^- s)^2}{\bar{M}(t-s)}} \\
&\times (1 + |y - a_k^- s|)^{-3/4} s^{-1/2} (1+s)^{-3/4} dy ds.
\end{aligned}
\end{equation*}
In the event that $|x| \ge |a_k^-| t$, we write
\begin{equation} \label{tminussdecomp}
x - a_j^- (t-s) -a_k^- s = (x - a_k^- t) - (a_j^- - a_k^-) (t-s),
\end{equation}
for which there is no cancellation between the summands, and we can
proceed exactly as in (\ref{cancellationargument}).  On the other
hand, for $|x| \le |a_j^-| t$, we write
\begin{equation} \label{sdecomp}
x - a_j^-(t-s) -a_k^- s = (x - a_j^- t) - (a_k^- - a_s^-) s,
\end{equation}
for which again there is no cancellation between summands, and we
can proceed exactly as in (\ref{cancellationargument}). For the
critical case $|a_j^-| t \le |x| \le |a_k^-| t$, we subdivide the
analysis into cases: $s \in [0,t/2]$ and $s \in [t/2,t]$. For $s \in
[0,t/2]$, we observe through (\ref{sdecomp}) the estimate
\begin{equation} \label{kernels}
\begin{aligned}
&e^{-\frac{(x - a_j^- (t-s) - a_k^- s)^2}{\bar{M}(t-s)}} (1+s)^{-3/4} \\
&\le C \Big[ e^{-\frac{(x - a_j^- t)^2}{\bar{M}(t-s)}} (1+s)^{-3/4}
+ e^{-\frac{(x - a_j^- (t-s) - a_k^- s)^2}{\bar{M}(t-s)}} (1 + |x -
a_j^- t|)^{-3/4} \Big].
\end{aligned}
\end{equation}
For the first, we can proceed exactly as in
(\ref{cancellationargument}), while for the second, we integrate
$(1+|y-a_k^- s|)^{-3/4}$ to obtain the estimate
\begin{equation*}
\begin{aligned}
C_1 &t^{-1} (1 + |x - a_j^- t|)^{-3/4} \int_0^{t/2}  e^{-\frac{(x -
a_j^- (t-s) - a_k^- s)^2}{\bar{M}(t-s)}}
s^{-1/2} (1+s)^{1/4} ds \\
&\le C t^{-1/2}  (1 + |x - a_j t|)^{-3/4},
\end{aligned}
\end{equation*}
where in establishing this last inequality we have integrated
\begin{equation*}
e^{-\frac{(x - a_j^- (t-s) - a_k^- s)^2}{\bar{M}(t-s)}}
\end{equation*}
in $s$. For $s \in [t/2, t]$, we observe through
(\ref{tminussdecomp}) the estimate
\begin{equation} \label{kerneltminuss}
\begin{aligned}
(t-s)^{-3/4} &e^{-\frac{(x - a_j^- (t-s) - a_k^- s)^2}{\bar{M}(t-s)}} \\
&\le C \Big[ (t-s)^{-3/4} e^{-\frac{(x - a_k^- t)^2}{\bar{M}(t-s)}}
+ |x - a_k^- t|^{-3/4} e^{-\frac{(x - a_j^- (t-s) - a_k^-
s)^2}{\bar{M}(t-s)}} \Big].
\end{aligned}
\end{equation}
For the first of these estimates, we proceed exactly as in
(\ref{cancellationargument}), while for the second, we integrate the
kernel in $y$ to obtain an estimate by
\begin{equation*}
\begin{aligned}
C_2 &|x - a_k^- t|^{-3/4} t^{-1/2} (1+t)^{-3/4}
\int_{t/2}^t (t-s)^{1/4} e^{-\frac{(x - a_j^- (t-s) - a_k^- s)^2}{\bar{M}(t-s)}} ds\\
&\le C |x - a_k^- t|^{-3/4} (1+t)^{-1/2}.
\end{aligned}
\end{equation*}
The apparent blowup as $x$ approaches $a_k^- t$ can be removed
similarly as in (\ref{smalltime}).

For the second estimate in (\ref{case1balance}), we have
\begin{equation*}
\begin{aligned}
\int_0^t &\int_{-|a_1^-|s}^0 (t-s)^{-1} e^{-\frac{(x - y - a_j^- (t-s))^2}{M(t-s)}} \\
& \times (1 + |y - a_k^- s| + |x - a_j^- (t-s) -  a_k^- s|)^{-3/4}
s^{-1/2} (1+s)^{-3/4} dy ds.
\end{aligned}
\end{equation*}
Again, we proceed in detail, only in Case (ii), $a_k^- < a_j^- < 0$.
In the event that $|x| \ge |a_k^-| t$ (though keeping in mind that
we have already considered the case $|x| \ge |a_1^-| t$), we use
(\ref{tminussdecomp}), for which there is no cancellation between
the summands, and upon integration of the kernel, we have an
estimate by
\begin{equation}
\begin{aligned} \label{nocanc}
C &(1+|x-a_k^- t|)^{-3/4} \int_0^t (t-s)^{-1/2}s^{-1/2} (1+s)^{-3/4} ds \\
&\le C_1 (1+|x-a_k^- t|)^{-3/4} t^{-1/2} \int_0^{t/2} s^{-1/2}(1+s)^{-3/4} ds \\
+ &C_2 (1+|x-a_k^- t|)^{-3/4}t^{-1/2} (1+t)^{-3/4} \int_{t/2}^t
(t-s)^{-1/2} ds,
\end{aligned}
\end{equation}
from which we obtain an estimate by $\bar\psi_1$.  Likewise, in the
case $|x| \le |a_j^-| t$, we use (\ref{sdecomp}), for which there is
no cancellation between summands, and we can proceed similarly as in
(\ref{nocanc}). For the critical case $|a_j^-| t \le |x| \le |a_k^-|
t$, we subdivide the analysis into cases: $s \in [0,t/2]$ and $s \in
[t/2,t]$. For $s \in [0,t/2]$, we observe through (\ref{sdecomp})
the estimate
\begin{equation*}
\begin{aligned}
&(1 + |y - a_k^- s| +  |x - a_j^- (t-s) - a_k^- s|)^{-3/4} (1 + s)^{-3/4} \\
& \le C \Big[(1 + |x - a_j^- t|)^{-3/4} (1 + s)^{-3/4} \\
&+ (1 + |y - a_k^- s| + |x - a_j^- (t-s) - a_k^- s|)^{-3/4} (1 + s +
|x - a_j^- t|)^{-3/4} \Big].
\end{aligned}
\end{equation*}
For the first of these estimates, we can proceed as in
(\ref{nocanc}), while for the second, we integrate the kernel in $y$
to obtain an estimate by
\begin{equation*}
\begin{aligned}
C_1 &(1+|x - a_j^- t|)^{-3/4} t^{-1/2} \int_0^{t/2}
(1 + |x - a_j^- (t-s) - a_k^- s|)^{-3/4} s^{-1/2} ds \\
&\le C (1+t)^{-1/2} (1+|x - a_j t|)^{-3/4}.
\end{aligned}
\end{equation*}
For $s \in [t/2, t]$, we observe through (\ref{tminussdecomp}) the
estimate
\begin{equation*}
\begin{aligned}
& (t-s)^{-3/4} (1 + |y - a_k^- s| +  |x - a_j^- (t-s) - a_k^- s|)^{-3/4}  \\
& \le (t-s)^{-3/4} \Big[(1 + |x - a_k^- t|)^{-3/4} \\
&+|x - a_k^- t|^{-3/4} (1 + |y - a_k^- s| + |x - a_j^- (t-s) - a_k^-
s|)^{-3/4} \Big].
\end{aligned}
\end{equation*}
For the first of these estimates, we proceed exactly as in
(\ref{cancellationargument}), while for the second, we integrate the
kernel in $y$ to obtain an estimate by
\begin{equation*}
\begin{aligned}
C_2 &|x - a_k^- t|^{-3/4} t^{-1/2} (1 + t)^{-3/4} \int_{t/2}^t
(t-s)^{1/4}
(1 + |x - a_j^- (t-s) - a_k^- s|)^{-3/4} \\
&\le C |x - a_k^- t|^{-3/4} (1+t)^{-3/4}.
\end{aligned}
\end{equation*}

{\it Nonlinearity} $(1+s)^{-1/4} s^{-1/2} \psi_2$.  We next consider
integrals of the form
\begin{equation} \label{case2integral}
\int_0^t \int_{-\infty}^0 (t-s)^{-1} e^{-\frac{(x - y - a_j^-
(t-s))^2}{M (t-s)}} (1+|y-a_k^- s|+s^{1/2})^{-3/2}s^{-1/2}
(1+s)^{-1/4} dy ds.
\end{equation}
In this case, we again observe (\ref{standarddecomposition}), from
which we obtain the estimate
\begin{equation} \label{case2balance}
\begin{aligned}
&e^{-\frac{(x - y - a_j^- (t-s))^2}{M (t-s)}} (1+|y-a_k^- s|+s^{1/2})^{-3/2} \\
&\le C \Big[ e^{-\frac{(x - a_j^- (t-s) - a_k^- s)^2}{\bar{M}
(t-s)}} e^{-\epsilon \frac{(x - y - a_j^- (t-s))^2}{M (t-s)}}
(1+|y-a_k^- s|+s^{1/2})^{-3/2} \\
&+ e^{-\frac{(x - y - a_j^- (t-s))^2}{M (t-s)}} (1+|y-a_k^-
s|+|x-a_j^- (t-s) -a_k^- s|+s^{1/2})^{-3/2} \Big].
\end{aligned}
\end{equation}
As in our analysis of the nonlinearity $(1+s)^{-1/4} s^{-1/2}
\bar\psi_1$, we focus on the case $a_k^- < a_j^- < 0$, compared to
which the other cases are either similar or more straightforward.
For the first estimate in (\ref{case2balance}), we have
\begin{equation*}
\begin{aligned}
\int_0^t &\int_{-\infty}^0 (t-s)^{-1} e^{-\epsilon \frac{(x - y -
a_j^- (t-s))^2}{M (t-s)}}
e^{-\frac{(x - a_j^- (t-s) - a_k^- s)^2}{\bar{M} (t-s)}} \\
&\times (1+|y-a_k^- s|+s^{1/2})^{-3/2} s^{-1/2} (1+s)^{-1/4} dy ds.
\end{aligned}
\end{equation*}
In the event that $|x| \ge |a_k^-| t$, we employ
(\ref{tminussdecomp}), for which there is no cancellation between
the summands.  In this case, we have an estimate by
\begin{equation} \label{nocanc2}
\begin{aligned}
C_1 t^{-1} &e^{-\frac{(x - a_k^- t)^2}{Lt}}
\int_0^{t/2} (1+s^{1/2})^{-1/2} s^{-1/2} (1+s)^{-1/4} ds \\
&+ C_2 (1+t^{1/2})^{-3/2} t^{-1/2} (1+t)^{-1/4} e^{-\frac{(x - a_k^-
t)^2}{Lt}}
\int_{t/2}^t (t-s)^{-1/2} ds \\
&\le C t^{-1} \ln(e+t) e^{-\frac{(x - a_k^- t)^2}{Lt}},
\end{aligned}
\end{equation}
where the seeming blowup as $t \to 0$ can be eliminated as in
(\ref{smalltime}). In the case $|x| \le |a_j^-| t$, we employ
(\ref{sdecomp}) for which again there is no cancellation between the
summands and we can proceed as in (\ref{nocanc2}).

We turn next to the critical case $|a_j^-| t \le |x| \le |a_k^-| t$,
for which we divide the analysis into cases $s \in [0,t/2]$ and $s
\in [t/2, t]$.  For $s \in [0,t/2]$, we observe through
(\ref{sdecomp}) the inequality
\begin{equation*}
\begin{aligned}
&e^{-\frac{(x - a_j^- (t-s) - a_k^- s)^2}{\bar{M} (t-s)}} s^{-1/2} (1+s)^{-1/4} \\
&\le C \Big[ e^{-\frac{(x - a_j^- t)^2}{L t}}s^{-1/2} (1+s)^{-1/4} \\
&\quad + e^{-\frac{(x - a_j^- (t-s) - a_k^- s)^2}{\bar{M} (t-s)}} |x
- a_j^- t|^{-1/2} (1 + |x - a_j^- t|)^{-1/4} \Big].
\end{aligned}
\end{equation*}
For the first of these last two estimates, we can proceed as in
(\ref{nocanc2}), while for the second, we have an estimate of the
form
\begin{equation*}
\begin{aligned}
C_1 t^{-1}  &|x - a_j^- t|^{-1/2}(1+|x-a_j^- t|)^{-1/4} \int_0^{t/2}
e^{-\frac{(x - a_j^- (t-s) - a_k^- s)^2}{\bar{M} (t-s)}} (1+s^{1/2})^{-1/2} ds \\
&\le C t^{-1/2} |x - a_j t|^{-1/2}  (1+|x-a_j t|)^{-1/4},
\end{aligned}
\end{equation*}
where as usual the singular behavior can be removed. For $s \in
[t/2,t]$, we observe through (\ref{tminussdecomp}) the inequality
\begin{equation*}
\begin{aligned}
&(t-s)^{-3/4} e^{-\frac{(x - a_j^- (t-s) - a_k^- s)^2}{\bar{M} (t-s)}} \\
&\le C \Big[ |x - a_k^- t|^{-3/4} e^{-\frac{(x - a_j^- (t-s) - a_k^-
s)^2}{\bar{M} (t-s)}} + (t-s)^{-3/4} e^{-\frac{(x - a_k^- t)^2}{L
t}}.
\end{aligned}
\end{equation*}
For the second of these last two estimates, we can proceed as in
(\ref{nocanc2}), while for the first we have an estimate of the form
\begin{equation*}
\begin{aligned}
C_1 &|x - a_k^- t|^{-3/4} t^{-1/2} (1+t)^{-1/4} (1+t^{1/2})^{-3/2}
\int_{t/2}^t  (t-s)^{1/4} e^{-\frac{(x - a_j^- (t-s) - a_k^- s)^2}{\bar{M} (t-s)}} ds \\
&\le C |x - a_k^- t|^{-3/4} (1+t)^{-3/4}.
\end{aligned}
\end{equation*}

For the second estimate in (\ref{case2balance}), we have
\begin{equation*}
\begin{aligned}
\int_0^t &\int_{-\infty}^0 (t-s)^{-1} e^{-\frac{(x - y - a_j^- (t-s))^2}{M (t-s)}} \\
&\times (1+|y-a_k^- s|+|x-a_j^- (t-s) -a_k^- s|+s^{1/2})^{-3/2}
s^{-1/2} (1+s)^{-1/4} dy ds.
\end{aligned}
\end{equation*}
Focusing on the case $a_k^- < a_j^- < 0$, we observe as before that
for $|x| \ge |a_k^-| t$, we have no cancellation between the
summands in (\ref{tminussdecomp}), and consequently there holds
\begin{equation*}
\begin{aligned}
&(1+|y-a_k^- s|+|x-a_j^- (t-s) -a_k^- s|+s^{1/2})^{-3/2} \\
&\le C (1+|y-a_k^- s|+|x-a_k^- t| + (t-s) + s^{1/2})^{-3/2},
\end{aligned}
\end{equation*}
from which, upon integration of the kernel, we have an estimate by
\begin{equation} \label{nocanc2a}
\begin{aligned}
C_1 &(1+|x-a_k^- t|+t^{1/2})^{-3/2}
\int_0^t (t-s)^{-1/2}s^{-1/2} (1+s)^{-1/4} \\
&\le C (1+t)^{-1/4} (1+|x-a_k^- t|+t^{1/2})^{-3/2}.
\end{aligned}
\end{equation}
In the event that $|x| \le |a_j^-| t$, we have no cancellation
between the summands in (\ref{sdecomp}), and consequently there
holds
\begin{equation*}
\begin{aligned}
&(1+|y-a_k^- s|+|x-a_j^- (t-s) -a_k^- s|+s^{1/2})^{-3/2} \\
&\le C (1+|y-a_k^- s|+|x-a_k^- t| + s)^{-3/2},
\end{aligned}
\end{equation*}
from which we have an estimate by
\begin{equation} \label{nocanc2b}
\begin{aligned}
C_1 &t^{-1} (1+|x - a_j^- t|)^{-1/2} \int_0^{t/2} s^{-1/2}(1+s)^{-1/4} ds \\
& +
C_2 t^{-1/2} (1+t)^{-1/4} (1 + |x - a_j^- t| + t)^{-3/2} \int_{t/2}^t (t-s)^{-1/2} ds \\
&\le C t^{-3/4} (1+|x - a_j^- t|)^{-1/2}.
\end{aligned}
\end{equation}
In the critical case $|a_j^-| t \le |x| \le |a_k^-| t$, we subdivide
the analysis further into cases $s \in [0,t/2]$ and $s \in [t/2,
t]$.  In the case $s \in [0,t/2]$, we observe through
(\ref{sdecomp}) the estimate
\begin{equation*}
\begin{aligned}
&(1+|y-a_k^- s|+|x-a_j^- (t-s) -a_k^- s|+s^{1/2})^{-3/2} s^{-1/2} (1+s)^{-1/4} \\
&\le C \Big[ (1+|y-a_k^- s|+|x-a_j^- (t-s) -a_k^- s|+|x-a_j^- t| +
s^{1/2})^{-3/2} \\
& \quad \times
s^{-1/2} (1+s)^{-1/4} \\
&+ (1+|y-a_k^- s|+|x-a_j^- (t-s) -a_k^- s|+s^{1/2})^{-3/2} \\
&\times (s+|x-a_j^- t|)^{-1/2} (1+s+|x-a_j^- t|)^{-1/4} \Big].
\end{aligned}
\end{equation*}
For the first of these last two estimates, we can proceed as in
(\ref{nocanc2b}), while for the second, we obtain an estimate by
\begin{equation*}
\begin{aligned}
C_1& t^{-1} |x-a_j^- t|^{-1/2} (1+|x-a_j^- t|)^{-1/4} \\
& \quad \times \int_0^{t/2} (1+|x-a_j^- (t-s) -a_k^- s|+s^{1/2})^{-1/2} ds \\
&\le C t^{-1/2} |x-a_j^- t|^{-1/2} (1+|x-a_j t|)^{-1/4},
\end{aligned}
\end{equation*}
where the apparent singularities can be eliminated as in
(\ref{smalltime}). In the case $s \in [t/2, t]$, we observe through
(\ref{tminussdecomp}) the estimate
\begin{equation*}
\begin{aligned}
&(t-s)^{-3/4} (1+|y-a_k^- s|+|x-a_j^- (t-s) -a_k^- s|+s^{1/2})^{-3/2} \\
&\le C \Big[|x - a_k^- t|^{-3/4} (1+|y-a_k^- s|+|x-a_j^- (t-s) -a_k^- s| + s^{1/2})^{-3/2} \\
&+ (t-s)^{-3/4} (1+|y-a_k^- s|+|x-a_k^- t + s^{1/2})^{-3/2} \Big].
\end{aligned}
\end{equation*}
For the second of the last two estimates, keeping in mind that here
$s \in [t/2, t]$, we can proceed as in (\ref{nocanc2a}), while for
the first we integrate the kernel to obtain an estimate by
\begin{equation*}
\begin{aligned}
C_1& |x-a_k^- t|^{-3/4} t^{-1/2} (1+t)^{-1/4} \\
&\times \int_{t/2}^t (t-s)^{1/4} (1+|x-a_j^- (t-s) -a_k^- s|+|x-a_j^- t| + s^{1/2})^{-3/2} ds \\
&\le C |x - a_k^- t|^{-3/4}t^{-1/2} (1+t)^{-1/4},
\end{aligned}
\end{equation*}
where the apparent singularities can be eliminated as in
(\ref{smalltime}).

{\it Nonlinearity} $(1+s)^{-1/4} s^{-1/2} \alpha$. We next consider
integrals of the form
\begin{equation} \label{case3integral}
\int_0^t \int_{-|a_1^-| s}^0 (t-s)^{-1} e^{-\frac{(x-y-a_j^-
(t-s))^2}{M(t-s)}} (1+s)^{-1} s^{-1/2} (1+|y|)^{-1/2} dy ds.
\end{equation}
In the case $|x| \ge |a_1^-| t$, we observe that there is no
cancellation between summands in decomposition (\ref{firstdecomp}),
and consequently, we have an estimate by
\begin{equation} \label{nocanc3}
\begin{aligned}
C_1 &t^{-1} e^{-\frac{(x-a_1^- t)^2}{Lt}} \int_0^{t/2} (1+s)^{-1/2} s^{-1/2} ds \\
&+ C_2 (1+t)^{-1} t^{-1/2} e^{-\frac{(x-a_1^- t)^2}{Lt}} \int_{t/2}^t (t-s)^{-1/2} ds \\
&\le C t^{-1} \ln (e+t) e^{-\frac{(x-a_1 t)^2}{Lt}},
\end{aligned}
\end{equation}
where the apparent blowup as $t \to 0$ can be removed as in
(\ref{smalltime}). For the remainder of the analysis, we restrict
our attention to the case $|x| \le |a_1^-| t$ and $a_j^- < 0$ (the
case $a_j^- > 0$ is similar and more direct).  Observing that for $y
\in [-|a_1^-| s,0]$, we have
\begin{equation*}
(1+s)^{-1} \le C (1 + s + |y|)^{-1},
\end{equation*}
we have the inequality
\begin{equation} \label{case3balance}
\begin{aligned}
&e^{-\frac{(x-y-a_j^- (t-s))^2}{M(t-s)}}
(1+s)^{-1} s^{-1/2} (1+|y|)^{-1/2} \Big{|}_{\{y \in [-|a_1^-| s,0]\}} \\
&\le C \Big[e^{-\epsilon \frac{(x-y-a_j^- (t-s))^2}{M(t-s)}}
e^{-\frac{(x-a_j^- (t-s))^2}{\bar{M}(t-s)}}
(1+s)^{-1} s^{-1/2} (1+|y|)^{-1/2} \\
&+ e^{-\frac{(x-y-a_j^- (t-s))^2}{M(t-s)}}
(1+s+|y|+|x-a_j^- (t-s)|)^{-1} \\
& \times (s+|y|+|x-a_j^- (t-s)|)^{-1/2} (1+|y|+|x-a_j^-
(t-s)|)^{-1/2} \Big].
\end{aligned}
\end{equation}
For the first estimate in (\ref{case3balance}), we observe that for
$|x| \ge |a_j^-| t$, there is no cancellation between $x-a_j^- t$
and $a_j^- s$, and consequently we have an estimate by
\begin{equation} \label{nocanc3a}
\begin{aligned}
C_1 &t^{-1} e^{-\frac{(x-a_j^- t)^2}{Lt}} \int_0^{t/2} (1+s)^{-1/2} s^{-1/2} ds \\
+ &C_2 (1+t)^{-1} t^{-1/2} e^{-\frac{(x-a_j^- t)^2}{Lt}}
\int_{t/2}^t (t-s)^{-1/2} ds \\
&\le C t^{-1} \ln(e+t) e^{-\frac{(x-a_j^- t)^2}{Lt}},
\end{aligned}
\end{equation}
where the apparent blowup as $t \to 0$ can be removed as in
(\ref{smalltime}). In the case $|x| \le |a_j^-| t$, we divide the
analysis into the cases, $s \in [0,t/2]$ and $s \in [t/2, t]$.  For
$s \in [0,t/2]$, we have the inequality
\begin{equation*}
\begin{aligned}
&e^{-\frac{(x-a_j^- (t-s))^2}{\bar{M}(t-s)}}
(1+s)^{-1} s^{-1/2} \\
&\le C\Big[ e^{-\frac{(x-a_j^- t)^2}{Lt}} e^{-\epsilon
\frac{(x-a_j^- (t-s))^2}{\bar{M}(t-s)}}
(1+s)^{-1} s^{-1/2} \\
&+ e^{-\frac{(x-a_j^- (t-s))^2}{\bar{M}(t-s)}} (1+s+|x-a_j^-
t|)^{-1} (s+|x-a_j^- t|)^{-1/2} \Big]
\end{aligned}
\end{equation*}
For the first of these last two estimates, we proceed as in
(\ref{nocanc3a}), while for the second, upon integration of
$(1+|y|)^{-1/2}$, we have an estimate on (\ref{case3integral}) by
\begin{equation*}
\begin{aligned}
C_1 &t^{-1} (1+|x - a_j^- t|)^{-3/4}
\int_0^{t/2} e^{-\frac{(x-a_j^- (t-s))^2}{\bar{M}(t-s)}} \\
&\times (1+s+|x-a_j^- t|)^{-1/4} (s+|x-a_j^- t|)^{-1/2} (1+s)^{1/2} ds \\
&\le Ct^{-1/2} (1+|x - a_j^- t|)^{-3/4}.
\end{aligned}
\end{equation*}
In the case $s \in [t/2, t]$, we observe the inequality
\begin{equation*}
(t-s)^{-1/2} e^{-\frac{(x-a_j^- (t-s))^2}{\bar{M}(t-s)}} \le C
|x|^{-1/2} e^{-\frac{(x-a_j^- (t-s))^2}{\bar{M}(t-s)}},
\end{equation*}
from which, upon integration of the kernel, we obtain an estimate on
(\ref{case3integral}) by
\begin{equation*}
\begin{aligned}
C_2 &(1+t)^{-1} t^{-1/2} |x|^{-1/2}
\int_{t/2}^t e^{-\frac{(x-a_j (t-s))^2}{\bar{M}(t-s)}} ds \\
&\le C t^{-1} |x|^{-1/2}.
\end{aligned}
\end{equation*}
For the second estimate in (\ref{case3balance}), we observe that in
the case $|x| \ge |a_j^-| t$, we have no cancellation between
$(x-a_j^- t)$ and $a_j^- s$, and consequently can estimate
\begin{equation} \label{nocanc3b}
\begin{aligned}
C_1 &t^{-1/2} (1+|x-a_j^- t|)^{-3/2}
\int_0^{t/2} s^{-1/2} ds \\
+ & C_2 t^{-1/2} (1 + |x - a_j^- t|)^{-3/2}
\int_{t/2}^t (t-s)^{-1/2} \\
&\le C (1+|x-a_j^- t|)^{-3/2},
\end{aligned}
\end{equation}
which suffices since the case $|x-a_j^- t| \le C\sqrt{t}$ requires
only $t^{-3/4}$ decay. In the case $|x| \le |a_j^-| t$, we subdivide
the analysis further into the cases $s \in [0,t/2]$ and $s \in
[t/2,t]$.  For $s \in [0, t/2]$, we observe the inequality
\begin{equation*}
(1+s+|x-a_j^- (t-s)|)^{-1} \le C (1+s+|x-a_j^- t|)^{-1},
\end{equation*}
from which, upon integration of the kernel, we obtain an estimate
by
\begin{equation*}
C_1 t^{-1/2} (1+|x-a_j^- t|)^{-1} \int_0^{t/2} s^{-1/2}
(1 + |x - a_j^- (t-s)|) ds \le C
t^{-1/2} (1+|x-a_j^- t|)^{-1}.
\end{equation*}
For $s \in [t/2, t]$, we observe the inequality
\begin{equation*}
\begin{aligned}
&(t-s)^{-1/2} (1+|y|+|x-a_j^- (t-s)|)^{-1/2} \\
&\le \Big[|x|^{-1/2} (1+|y|+|x-a_j^- (t-s)|)^{-1/2} + (t-s)^{-1/2}
(1+|y|+|x|)^{-1/2}\Big].
\end{aligned}
\end{equation*}
For the first of these last two estimates, we obtain an estimate on
(\ref{case3integral}) by
\begin{equation*}
\begin{aligned}
C_2 (1+t)^{-1} t^{-1/2} |x|^{-1/2} \int_{t/2}^t (1+|x-a_j^-
(t-s)|)^{-1/2} ds \le C (1+t)^{-1} |x|^{-1/2},
\end{aligned}
\end{equation*}
while for the second, we have an estimate by
\begin{equation*}
\begin{aligned}
C_2 (1+t)^{-1} t^{-1/2} (1+|x|)^{-1/2} \int_{t/2}^t (t-s)^{-1/2} ds
\le C (1+t)^{-1} |x|^{-1/2}.
\end{aligned}
\end{equation*}

{\it Nonlinearity} $s^{-1/2} (1+s)^{-1/4} \varphi$.  We next
consider integrals of the form
\begin{equation} \label{case4integral}
\int_0^t \int_{-\infty}^0 (t-s)^{-1} e^{-\frac{(x - y - a_j^-
(t-s))^2}{M (t-s)}} (1+s)^{-3/4} s^{-1/2} e^{-\frac{(y - a_k^-
s)^2}{Ms}} dy ds.
\end{equation}
In this case we observe from Lemma 6 of \cite{HZ} the equality
\begin{equation} \label{completedsquare}
\begin{aligned}
& e^{-\frac{(x - y - a_j^- (t-s))^2}{M (t-s)}} e^{-\frac{(y - a_k^- s)^2}{Ms}} \\
&= e^{-\frac{(x-a_j^- (t-s)-a_k^- s)^2}{Mt}} e^{-\frac{t}{Ms(t-s)}(y
- \frac{xs-(a_j^-+a_k^-)(t-s)s}{t})^2},
\end{aligned}
\end{equation}
from which direct integration over $y$ leads to an estimate by
\begin{equation*}
C t^{-1/2} \int_0^t (t-s)^{-1/2} (1+s)^{-3/4} e^{-\frac{(x-a_j^-
(t-s)-a_k^- s)^2}{Mt}}.
\end{equation*}
As in the previous analyses, we focus on the case $a_k^- < a_j^-
<0$, and note that analysis of the remaining case $a_j^- \le a_k^- <
0$ is similar.  In the event that $|x| \ge |a_k^-| t$, we observe
that there is no cancellation between the summands of
(\ref{tminussdecomp}) and we have an estimate by
\begin{equation} \label{nocanc4}
\begin{aligned}
C_1& t^{-1} e^{-\frac{(x-a_k^- t)^2}{Lt}} \int_0^{t/2} (1+s)^{-3/4} ds \\
&+ C_2 t^{-1/2} (1+t)^{-3/4} e^{-\frac{(x-a_k^- t)^2}{Lt}}
\int_{t/2}^t (t-s)^{-1/2} ds \\
&\le C t^{-3/4} e^{-\frac{(x-a_k^- t)^2}{Lt}}.
\end{aligned}
\end{equation}
For $|x| \le |a_k^-| t$, we divide the analysis into cases, $s \in
[0,t/2]$ and $s \in [t/2, t]$.  For $s \in [0,t/2]$, use
(\ref{sdecomp}), for which the summands cancel and we have the
estimate (\ref{kernels}).  For the first estimate in
(\ref{kernels}), we can proceed exactly as in (\ref{nocanc4}), while
for the second, we have an estimate by
\begin{equation*}
\begin{aligned}
C_1& t^{-1} (1+|x-a_j^- t|)^{-3/4}
\int_0^{t/2} e^{-\frac{(x-a_j^- (t-s)-a_k^- s)^2}{Mt}} ds \\
&\le C t^{-1/2} (1+|x-a_j t|)^{-3/4},
\end{aligned}
\end{equation*}
where the apparent blowup as $t \to 0$ can be removed as in
(\ref{smalltime}). For $s \in [t/2, t]$, we observe through
(\ref{tminussdecomp}) the estimate (\ref{kerneltminuss}). For the
first estimate in (\ref{kerneltminuss}), we proceed as in
(\ref{nocanc4}), while for the second we have an estimate by
\begin{equation*}
\begin{aligned}
C_2 &t^{-1/2} |x - a_k^- t|^{-3/4} (1+t)^{-3/4}
\int_{t/2}^t (t-s)^{1/4} e^{-\frac{(x-a_j^- (t-s)-a_k^- s)^2}{Mt}} ds \\
&\le C t^{-1/2} |x - a_k^- t|^{-3/4}.
\end{aligned}
\end{equation*}

{\it Nonlinearity} $(1+s)^{-1} e^{-\eta |y|}$.  We next consider
integrals of the form
\begin{equation} \label{case5integral}
\int_0^t \int_{-\infty}^0 (t-s)^{-1} e^{-\frac{(x - y - a_j^-
(t-s))^2}{M (t-s)}} e^{-\eta|y|} (1+s)^{-1} dy ds.
\end{equation}
First, we observe that in the event $|y| \ge |a_1^-| s$, we have
exponential decay in both $y$ and $s$, from which the claimed
estimate readily follows.  For what remains, then, we focus on the
case $|y| \le |a_1^-| s$. In the event that $|x| \ge |a_1^-| t$, we
write
\begin{equation*}
x - y - a_j^- (t-s) = (x - a_1^- t) - (y - a_1^- s) - (a_j^- -
a_1^-) (t-s),
\end{equation*}
for which there is no cancellation between the summands (in the case
$y \in [a_1^- s, 0]$), and we immediately arrive at an estimate by
\begin{equation*}
C t^{-1} e^{-\frac{(x - a_1^- t)^2}{Lt}}.
\end{equation*}
In the case $|x| \le |a_1^-| t$, we focus on the critical case
$a_j^- < 0$.  Here, we observe the estimate
\begin{equation}
\begin{aligned} \label{case5balance}
&e^{-\frac{(x-y-a_j^- (t-s))^2}{M(t-s)}}e^{-\eta |y|} \\
&\le C \Big[e^{-\frac{(x-a_j^- (t-s))^2}{\bar{M}(t-s)}}e^{-\eta |y|}
+ e^{-\eta_1 |x-a_j^- (t-s)|} e^{-\eta_2 |y|}\Big],
\end{aligned}
\end{equation}
for some constants $\bar{M}>0$, $\eta_1 > 0$, and $\eta_2 > 0$. For
the first estimate in (\ref{case5balance}), we have
\begin{equation*}
\int_0^t \int_{-\infty}^0 (t-s)^{-1} e^{-\frac{(x-a_j^-
(t-s))^2}{\bar{M}(t-s)}}e^{-\eta |y|} (1+s)^{-1} dy ds.
\end{equation*}
In the event that $|x| \ge |a_j^-| t$, we write
\begin{equation} \label{decomp5a}
x - a_j^- (t-s) = (x - a_j^- t) + a_j^- s,
\end{equation}
for which there is no cancellation between the summands, and we
arrive at an estimate by
\begin{equation} \label{nocanc5}
\begin{aligned}
C_1& t^{-1} e^{-\frac{(x - a_j^- t)^2}{Lt}} \int_0^{t/2} (1+s)^{-1}
ds
+ C_2 (1+t)^{-1} e^{-\frac{(x - a_j^- t)^2}{Lt}} \int_{t/2}^{t-1} (t-s)^{-1} ds \\
&+ C_3 (1+t)^{-1} e^{-\frac{(x - a_j^- t)^2}{Lt}} \int_{t-1}^{t}
(t-s)^{-1/2} ds \le C t^{-1} \ln(e+t) e^{-\frac{(x - a_j^-
t)^2}{Lt}}.
\end{aligned}
\end{equation}
For $|x| \le |a_j^-| t$, we divide the analysis into cases, $s \in
[0,t/2]$ and $s \in [t/2,t]$.  In the case $s \in [0,t/2]$, we
observe the cancellation between summands in (\ref{decomp5a}), which
leads to the estimate
\begin{equation*}
\begin{aligned}
&e^{-\frac{(x - a_j^- (t-s))^2}{\bar{M}(t-s)}} (1+s)^{-1} \\
&\le C\Big[ e^{-\frac{(x - a_j^- t)^2}{L t}} (1+s)^{-1} +
e^{-\frac{(x - a_j^- (t-s))^2}{\bar{M}(t-s)}} (1+|x - a_j^-
t|)^{-1}\Big].
\end{aligned}
\end{equation*}
For the first of these last two estimates, we proceed exactly as in
(\ref{nocanc5}), while for the second, we have an estimate by
\begin{equation*}
C_1 t^{-1} (1+|x - a_j^- t|)^{-1} \int_0^{t/2} e^{-\frac{(x - a_j^-
(t-s))^2}{\bar{M}(t-s)}} ds \le C t^{-1/2} (1+|x - a_j^- t|)^{-1}.
\end{equation*}
For $s \in [t/2, t]$, we observe that for $|a_j^-| (t-s) \le (1/2)
|x|$, we have the estimate
\begin{equation*}
(t-s)^{-1/2} e^{-\frac{(x-a_j^- (t-s))^2}{\bar{M} (t-s)}} \le C
(t-s)^{-1/2} e^{-\frac{x^2}{L |x|}} e^{-\frac{(x-a_j^-
(t-s))^2}{\bar{M} (t-s)}},
\end{equation*}
while for $|a_j^-| (t-s) \ge (1/2) |x|$, we have the estimate
\begin{equation*}
(t-s)^{-1/2} e^{-\frac{(x-a_j^- (t-s))^2}{\bar{M} (t-s)}} \le C
|x|^{-1/2} e^{-\frac{(x-a_j^- (t-s))^2}{\bar{M} (t-s)}}.
\end{equation*}
For the second of these last two estimates, we obtain an estimate
\begin{equation*}
C_2 (1+t)^{-1} |x|^{-1/2} \int_{t/2}^t (t-s)^{-1/2}
e^{-\frac{(x-a_j^- (t-s))^2}{\bar{M} (t-s)}} ds \le C (1+t)^{-3/4}
|x|^{-1/2}.
\end{equation*}
For the first, we have the same decay in $t$, with exponential decay
in $|x|$. For the second estimate in (\ref{case5balance}), we have
\begin{equation*}
\int_0^t \int_{-\infty}^0 (t-s)^{-1} e^{-\eta_1 |x-a_j^- (t-s)|}
e^{-\eta_2 |y|} (1+s)^{-1} dy ds.
\end{equation*}
In the event that $|x| \ge |a_j^-| t$, we observe that there is no
cancellation between the summands of (\ref{decomp5a}), so that
\begin{equation*}
e^{-\eta_1 |x-a_j^- (t-s)|} \le C e^{-\eta_3 |x-a_j^- t|} e^{-\eta_4
|a_j^-| s},
\end{equation*}
which leads immediately to an estimate better than $\psi_2$. For
$|x| \le |a_j^-| t$, we proceed similarly as with the first estimate
in (\ref{case5balance}).

This concludes our analysis of the main case, the first estimate in
Lemma \ref{nonlinearestimates}.

{\it Second estimate of Lemma \ref{nonlinearestimates}}.  For the
second estimate in Lemma \ref{nonlinearestimates}, we consider in
detail the case of nonlinearity $\partial_y (\varphi^i)^2$, for
which the other cases are similar.  Following the analysis of
\cite{Ra}, we divide the integration over $s$ as
\begin{equation} \label{diffusionwaves1}
\begin{aligned}
\int_0^t& \int_{-\infty}^{+\infty} \tilde{G} (x,t-s;y) \partial_y
(\varphi^k)^2 dy ds = \int_0^{\sqrt{t}} \int_{-\infty}^{+\infty}
\tilde{G} (x,t-s;y) \partial_y (\varphi^k)^2 dy ds \\
&+ \int_{\sqrt{t}}^{t-\sqrt{t}} \int_{-\infty}^{+\infty}
\tilde{G} (x,t-s;y) \partial_y (\varphi^k)^2 dy ds \\
& \quad + \int_{t-\sqrt{t}}^t \int_{-\infty}^{+\infty} \tilde{G}
(x,t-s;y) \partial_y (\varphi^k)^2 dy ds.
\end{aligned}
\end{equation}
For the first integral in (\ref{diffusionwaves1}), we focus on the
case $y < 0$ and on the the convection term of the Green's function,
for which we must estimate integrals of the form
\begin{equation*}
\int_0^{\sqrt{t}} \int_{-\infty}^0 (t-s)^{-1} e^{-\frac{(x-y-a_j^-
(t-s))^2}{M (t-s)}} (1+s)^{-1} e^{-\frac{(y-a_k^- s)^2}{Ms}} dy ds.
\end{equation*}
According to Lemma 6 from \cite{HZ}, we can write
(\ref{completedsquare}), from which direct integration over $y$
leads to an estimate by
\begin{equation} \label{diffusionwaves2}
Ct^{-1/2} \int_0^{\sqrt{t}} (t-s)^{-1/2} (1+s)^{-1/2} e^{-\frac{(x -
a_j^- (t-s) - a_k^- s)^2}{Mt}} ds.
\end{equation}
We have three cases to consider: (1) $a_k^- < 0 < a_j^-$, (2) $a_k^-
< a_j^- < 0$, and (3) $a_j^- < a_k^- < 0$; here we point out that
the case $a_k^- = a_j^-$ does not arise: the reason is that,
according to (\ref{ajabaa}), the $(\fe^k)^2$ term in $\Phi$ occurs
only in $r^-_j$ direction for $j\ne k.$ Hence if we consider the
detailed
 description of $S$ in Green function bounds (\ref{Sov}), we can see that the case
 $a_k^- = a_j^-$ does not arise (since $l_j^- r_k^- = 0$).   We will focus on the case $a_k^- < a_j^- < 0$, for which
the remaining cases are similar. First, in the event that $|x -
a_j^- t| \le C\sqrt{t}$ we can conclude decay of the form
\begin{equation*}
e^{-\frac{(x-a_j^- t)^2}{Mt}}
\end{equation*}
by boundedness. In the case $|x - a_j^- t| \ge C\sqrt{t}$, for $C$
sufficiently large, and for $s \in [0, \sqrt{t}]$, we have
\begin{equation*}
|x - a_j^- (t-s) - a_k^- s| = |(x - a_j^- t) - (a_k^- - a_j^-) s|
\ge c_1 |x-a_j^- t|,
\end{equation*}
for which we can estimate (\ref{diffusionwaves2}) for $t > 1$ by
\begin{equation} \label{nocancdiffusionwaves1}
C_1 t^{-1} e^{-\frac{(x-a_j^- t)^2}{Lt}} \int_0^{\sqrt{t}}
(1+s)^{-1/2} \le C t^{-1} (1+\sqrt{t})^{1/2} e^{-\frac{(x-a_j^-
t)^2}{Lt}}.
\end{equation}
For the third integral in (\ref{diffusionwaves1}), we can proceed
similarly to obtain an estimate by
\begin{equation*}
C t^{-1/4} (1+t)^{-1/2} e^{-\frac{(x-a_k^- t)^2}{Lt}}.
\end{equation*}
In the case of the second integral in (\ref{diffusionwaves1}), we
closely follow the approach of \cite{Liu}, using the framework
of \cite{Ra}.  Defining for some fixed
$k$ $\phi (x,t)$ as a Burgers kernel, we have
\begin{equation*}
\begin{cases}
\phi_t - \beta_k^- \phi_{xx} = -\gamma_k^- (\phi^2)_x& \text{for } t > -1, \\
\phi (x + a_k^-, 1) = m_k \delta_0& t=-1,
\end{cases}
\end{equation*}
for which we have $\varphi^k (x,t) = \phi (x - a_k^- t, t)$, where
the $\varphi^k (x, t)$ are as in (\ref{burgers}).  In this notation,
the second integral in (\ref{diffusionwaves1}) becomes
\begin{equation} \label{diffmain}
\int_{\sqrt{t}}^{t-\sqrt{t}} \int_{-\infty}^{+\infty} \tilde{G} (x,
t-s; y) (\phi(y-a_k^- s, s)^2)_y dy ds.
\end{equation}
The idea of Liu is to take advantage of the observation that
time derivatives of the heat kernel decay at the same algebraic rate
as two space derivatives, $t^{-3/2}$, and also of the relation
between time and space derivatives of $\tilde{G}$ and $\phi$.
Intuitively, we can think that we would like to integrate
(\ref{diffmain}) by parts in $s$ to put time derivatives on
$\tilde{G}$, but in order to facilitate this, we would like the
derivatives on $\phi^2$ to be with respect to $s$. Though convecting
kernels such as $\tilde{G}$ do not enjoy this property of
fast-decaying time derivatives, we observe that under a change of
integration variable, the leading order term in $\tilde{G}$ can be
converted into a heat kernel.  More precisely, according to
Proposition \ref{greenbounds}, the leading order contribution to the
scattering piece of the Green's function is given by $g (x - a_j^-
t, t)$, where $g (x, t)$ is the heat kernel
\begin{equation*}
g (x, t; y) := c t^{-1/2} e^{-\frac{(x - y)^2}{4\beta_j t}},
\end{equation*}
for a constant $c = r_j^- (l_j^-)^\text{tr}/\sqrt{4 \pi \beta_j^-}$.
Under the change of variables $\xi = y + a_j^- (t-s)$, the second
integral in (\ref{diffusionwaves1}) becomes, for this leading order
term,
\begin{equation*}
\int_{\sqrt{t}}^{t-\sqrt{t}} \int_{-\infty}^{+\infty} g (x, t-s;
\xi) (\phi(\xi - a_j^- (t-s) - a_k^- s, s)^2)_\xi d\xi ds.
\end{equation*}
Following \cite{Ra}, we write for $g = g(x,\tau,\xi)$ and $\phi =
\phi (\xi, \tau)$,
\begin{equation*}
\begin{aligned}
&\Big(g(x, t-s; \xi) (\phi(\xi - a_j^- (t-s) - a_k^- s, s)^2) \Big)_s \\
&=
- g_\tau (x, t-s; \xi) (\phi(\xi - a_j^- (t-s) - a_k^- s, s)^2) \\
&\quad + (a_j^- - a_k^-) g(x, t-s; \xi)
(\phi(\xi - a_j^- (t-s) - a_k^- s, s)^2)_\xi \\
&\quad \quad + g(x, t-s; \xi) (\phi(\xi - a_j^- (t-s) - a_k^- s,
s)^2)_\tau,
\end{aligned}
\end{equation*}
from which we have
\begin{equation} \label{Raoofidecomp}
\begin{aligned}
(a_j^- - a_k^-) &g(x, t-s; \xi) (\phi(\xi - a_j^- (t-s) - a_k^- s)^2)_\xi \\
=
&\Big(g(x, t-s; \xi) (\phi(\xi - a_j^- (t-s) - a_k^- s, s)^2) \Big)_s \\
&\quad +
g_\tau (x, t-s; \xi) (\phi(\xi - a_j^- (t-s) - a_k^- s, s)^2) \\
&-g(x, t-s; \xi) (\phi(\xi - a_j^- (t-s) - a_k^- s, s)^2)_\tau.
\end{aligned}
\end{equation}
We proceed now by analyzing (under integration) each term on the
right hand side of (\ref{Raoofidecomp}).  For the first, we have
\begin{equation} \label{Raoofidecomp1}
\begin{aligned}
&\int_{\sqrt{t}}^{t-\sqrt{t}} \int_{-\infty}^{+\infty}
\Big(g (x, t-s; \xi) (\phi(\xi - a_j^- (t-s) - a_k^- s, s)^2)\Big)_s d\xi ds \\
&=\int_{-\infty}^{+\infty}
g (x, \sqrt{t}; \xi) (\phi(\xi - a_j^- \sqrt{t} - a_k^- (t-\sqrt{t}), t-\sqrt{t})^2) d\xi \\
&+\int_{-\infty}^{+\infty} g (x, t-\sqrt{t}; \xi) (\phi(\xi - a_j^-
(t-\sqrt{t}) - a_k^- \sqrt{t}, \sqrt{t})^2) d\xi.
\end{aligned}
\end{equation}
For the first integral in (\ref{Raoofidecomp1}), we have an estimate
by
\begin{equation*}
C \int_{-\infty}^{+\infty} (\sqrt{t})^{-1/2} e^{-\frac{(x -
\xi)^2}{M \sqrt{t}}} (1+(t-\sqrt{t}))^{-1} e^{-\frac{(\xi - a_j^-
\sqrt{t} - a_k^- (t-\sqrt{t}))^2}{M (t-\sqrt{t})}} d\xi.
\end{equation*}
Writing now $\xi = z + a_j^- \sqrt{t}$, we can rewrite this last
integral as
\begin{equation*}
C \int_{-\infty}^{+\infty} (\sqrt{t})^{-1/2} e^{-\frac{(x - z -
a_j^- \sqrt{t})^2}{M \sqrt{t}}} (1+(t-\sqrt{t}))^{-1} e^{-\frac{(z
- a_k^- (t-\sqrt{t}))^2}{M (t-\sqrt{t})}} dz.
\end{equation*}
We have, then, according to Lemma 6 of \cite{HZ}, with $t-s$
replaced by $\sqrt{t}$,
\begin{equation*}
\begin{aligned}
&e^{-\frac{(x - z - a_j^- \sqrt{t})^2}{M \sqrt{t}}}
e^{-\frac{(z  - a_k^- (t-\sqrt{t}))^2}{M (t-\sqrt{t})}} \\
&= e^{-\frac{(x - a_j^- \sqrt{t} - a_k^- (t - \sqrt{t}))^2}{Mt}}
e^{-\frac{t}{M\sqrt{t} (t-\sqrt{t})} (\xi - \frac{xs-(a_j^- +
a_k^-)\sqrt{t}(t-\sqrt{t})^2}{t})}.
\end{aligned}
\end{equation*}
Upon integration in $\xi$, then, we have an estimate by
\begin{equation*}
C t^{-1/2} (1+(t-\sqrt{t}))^{-1} (t-\sqrt{t})^{1/2} e^{-\frac{(x -
a_j^- \sqrt{t} - a_k^- (t - \sqrt{t}))^2}{Mt}}.
\end{equation*}
In the event that $|x - a_k^- t| \le C \sqrt{t}$, we have
exponential decay
\begin{equation*}
e^{-\frac{(x - a_k^- t)^2}{Mt}}
\end{equation*}
from boundedness, while for $|x - a_k^- t| \ge C \sqrt{t}$, $C$
sufficiently large, we have
\begin{equation*}
|x - a_j^- \sqrt{t} - a_k^- (t - \sqrt{t})| = |(x - a_k^- t) -
(a_j^- - a_k^-) \sqrt{t}| \ge c |x - a_k^- t|.
\end{equation*}
In either case, we conclude for $t$ sufficiently large an estimate
\begin{equation*}
Ct^{-1} e^{-\frac{(x - a_k^- t)^2}{Lt}}.
\end{equation*}
The second integral in (\ref{Raoofidecomp1}) can be analyzed
similarly, and we obtain an estimate of the same form. For the
second integral arising from (\ref{Raoofidecomp}), we have
\begin{equation*}
\int_{\sqrt{t}}^{t-\sqrt{t}} \int_{-\infty}^{+\infty} (t-s)^{-3/2}
e^{-\frac{(x - \xi)^2}{M (t-s)}} (1+s)^{-1} e^{-\frac{(\xi - a_j^-
(t-s) - a_k^- s)^2}{Ms}} d\xi ds.
\end{equation*}
Returning to our original coordinates $\xi = y + a_j^- (t-s)$, we
have
\begin{equation*}
\int_{\sqrt{t}}^{t-\sqrt{t}} \int_{-\infty}^{+\infty} (t-s)^{-3/2}
e^{-\frac{(x - y - a_j^- (t-s))^2}{M (t-s)}} (1+s)^{-1} e^{-\frac{(y
- a_k^- s)^2}{Ms}} d\xi ds.
\end{equation*}
Proceeding now exactly as in the derivation of
(\ref{diffusionwaves2}), we arrive at an estimate by
\begin{equation} \label{diffusionwaves3}
C t^{-1/2} \int_{\sqrt{t}}^{t-\sqrt{t}} (t-s)^{-1} (1+s)^{-1/2}
e^{-\frac{(x - a_j^- (t-s) - a_k^- s)^2}{Mt}} ds.
\end{equation}
As in previous arguments, we will focus on the case $a_k^- < a_j^- <
0$, for which the remaining cases are similar.  First, in the event
that $|x| \ge |a_k^-| t$, we observe that there is no cancellation
between summands in (\ref{tminussdecomp}), and consequently that we
can estimate (\ref{diffusionwaves3}) for $t$ large enough so that
$\sqrt{t} \le t/2$ by
\begin{equation} \label{nocancdiffusionwaves3}
\begin{aligned}
C_1 &t^{-3/2} e^{-\frac{(x-a_k^- t)^2}{Lt}} \int_{\sqrt{t}}^{t/2} (1+s)^{-1/2} ds \\
&+ C_2 t^{-1/2} (1+t)^{-1/2} e^{-\frac{(x-a_k^- t)^2}{Lt}}
\int_{t/2}^{t-\sqrt{t}} (t-s)^{-1} e^{-\frac{(x - a_j^- (t-s) - a_k^- s)^2}{Mt}} ds  \\
&\le C t^{-1} e^{-\frac{(x-a_k^- t)^2}{Lt}}.
\end{aligned}
\end{equation}
For $|x| \le |a_j^-| t$, we observe that there is no cancellation
between summands in (\ref{sdecomp}), and consequently we can
estimate (\ref{diffusionwaves3}) similarly as in
(\ref{nocancdiffusionwaves3}) to obtain an estimate by
\begin{equation*}
C t^{-1} e^{-\frac{(x-a_j^- t)^2}{Lt}}.
\end{equation*}
For the critical case $|a_j^-| t \le |x| \le |a_k^-| t$, we divide
the analysis into the cases $s \in [\sqrt{t},t/2]$ and $s \in [t/2,
t-\sqrt{t}]$. For $s \in [\sqrt{t}, t/2]$, we observe through
(\ref{sdecomp}) the inequality
\begin{equation*}
\begin{aligned}
&e^{-\frac{(x - a_j^- (t-s) -a_k^- s)^2}{Mt}} (1+s)^{-1/2} \\
&\le C \Big[ e^{-\frac{(x-a_j^- t)^2}{Lt}} (1+s)^{-1/2} +
e^{-\frac{(x - a_j^- (t-s) -a_k^- s)^2}{Mt}} (1+|x-a_j^- t|)^{-1/2}
\Big].
\end{aligned}
\end{equation*}
For the first of these last two estimates, we proceed as in
(\ref{nocancdiffusionwaves3}), while for the second, we have an
estimate by
\begin{equation*}
C_1 t^{-3/2} (1+|x-a_j^- t|)^{-1/2} \int_{\sqrt{t}}^{t/2}
e^{-\frac{(x - a_j^- (t-s) -a_k^- s)^2}{Mt}} ds
\le C t^{-3/2} (1+t)^{1/2} (1+|x-a_j^- t|)^{-1/2}.
\end{equation*}
For $s \in [t/2, t-\sqrt{t}]$, we observe through
(\ref{tminussdecomp}) the estimate
\begin{equation*}
\begin{aligned}
(t-s)^{-1/2} &e^{-\frac{(x-a_j^- (t-s) - a_k^- s)^2}{Mt}} \\
&\le C \Big[ |x-a_k^- t|^{-1/2} e^{-\frac{(x-a_j^- (t-s) - a_k^-
s)^2}{Mt}} + (t-s)^{-1/2} e^{-\frac{(x - a_k^- t)^2}{Lt}} \Big].
\end{aligned}
\end{equation*}
For the second of these last two estimates we can proceed as in
(\ref{nocancdiffusionwaves3}) to determine an estimate by
\begin{equation*}
Ct^{-1} e^{-\frac{(x - a_k^- t)^2}{Lt}},
\end{equation*}
while for the first, we have an estimate by
\begin{equation*}
\begin{aligned}
C_1& t^{-1/2} (1+t)^{-1/2} |x - a_k^- t|^{-1/2}
\int_{\sqrt{t}}^{t-\sqrt{t}} (t-s)^{-1/2}
e^{-\frac{(x-a_j^- (t-s) - a_k^- s)^2}{Mt}} ds \\
&\le C (1+t)^{-1/2} |x - a_k^- t|^{-1/2} (\sqrt{t})^{-1/2}.
\end{aligned}
\end{equation*}
For the final expression on the right hand side of
(\ref{Raoofidecomp}), we have integrals of the form
\begin{equation*}
\int_{\sqrt{t}}^{t-\sqrt{t}} \int_{-\infty}^{+\infty} (t-s)^{-1/2}
e^{-\frac{(x - \xi)^2}{M (t-s)}} (1+s)^{-2} e^{-\frac{(\xi - a_j^-
(t-s) - a_k^- s)^2}{Ms}} d\xi ds.
\end{equation*}
Re-writing in our original variable $\xi = y + a_j^- (t-s)$, we have
\begin{equation*}
\int_{\sqrt{t}}^{t-\sqrt{t}} \int_{-\infty}^{+\infty} (t-s)^{-1/2}
e^{-\frac{(x - y - a_j^- (t-s))^2}{M (t-s)}} (1+s)^{-2} e^{-\frac{(y
- a_k^- s)^2}{Ms}} d\xi ds.
\end{equation*}
Proceeding now exactly as in the derivation of
(\ref{diffusionwaves3}), we arrive at an estimate by
\begin{equation} \label{diffusionwaves4}
C t^{-1/2} \int_{\sqrt{t}}^{t-\sqrt{t}} (1+s)^{-3/2} e^{-\frac{(x -
a_j^- (t-s) - a_k^- s)^2}{Mt}} ds.
\end{equation}
For the cases $|x| \ge |a_k^-| t$ and $|x| \le |a_j^-| t$, we can
proceed as in (\ref{nocancdiffusionwaves3}) to arrive at estimates
of the form $t^{-1/4} \varphi$, while in the case $|a_j^-| t \le |x|
\le |a_k^-| t$, we divide the analysis into case, $s \in
[\sqrt{t},t/2]$ and $s \in [t/2, t-\sqrt{t}]$. For $s \in [\sqrt{t},
t/2]$, we observe through (\ref{sdecomp}) the estimate
\begin{equation*}
\begin{aligned}
&e^{-\frac{(x - a_j^- (t-s) - a_k^- s)^2}{Mt}} (1+s)^{-3/2} \\
&\le C \Big[e^{-\frac{(x - a_j^- t)^2}{Lt}} (1+s)^{-3/2} +
(1+|x-a_j^- t| + t^{1/2})^{-3/2} e^{-\frac{(x - a_j^- (t-s) - a_k^-
s)^2}{Mt}} \Big].
\end{aligned}
\end{equation*}
For the first of these last two estimates, we proceed as in
(\ref{nocancdiffusionwaves3}), while for the second we have an
estimate by
\begin{equation*}
\begin{aligned}
C& t^{-1/2} (1 + |x-a_j^- t| + t^{1/2})^{-3/2}
\int_{\sqrt{t}}^{t/2} e^{-\frac{(x - a_j^- (t-s) - a_k^- s)^2}{Mt}} ds \\
&\le C (1 + |x-a_j^- t| + t^{1/2})^{-3/2}.
\end{aligned}
\end{equation*}
For the case $s \in [t/2, t-\sqrt{t}]$, we similarly arrive at an
estimate by $(1+t)^{-1} |x - a_k^- t|^{-1/2}$.

{\it Excited estimates.} We turn now to the estimates in Lemma
(\ref{nonlinearestimates}) that involve the {\it excited} terms $e_i
(y,t)$.  Observing that in the Lax and overcompressive cases, we
have the estimate
\begin{equation*}
|\partial_y e_i (y,t)| \le C t^{-1/2} \sum_{a_k^- > 0}
e^{-\frac{(y+a_k^- t)^2}{Mt}},
\end{equation*}
integrals in the third estimate of Lemma \ref{nonlinearestimates}
take the form
\begin{equation*}
\begin{aligned}
\int_0^t \int_{-|a_1^-| s}^0 (t-s)^{-1/2} e^{-\frac{(y+a_k^-
(t-s))^2}{M(t-s)}} \Psi (s,y) dy ds,
\end{aligned}
\end{equation*}
which differ from those of the previous analysis only by factor
$(t-s)^{-1/2}$. Proceeding almost exactly as in our analysis of the
first estimate of Lemma \ref{nonlinearestimates}, we determine an
estimate on these integrals of $C (1+t)^{-3/4}$.

{\it Excited diffusion wave estimates.}  The third and fourth
estimates in Lemma \ref{nonlinearestimates} are similar and
straightforward to prove, and we consider only the third.  As each
summand in $\Phi (y, s)$ can be dealt with similarly (see
(\ref{ghati})--(\ref{ajabaa})), we focus on the nonlinearity
$\partial_y (\varphi^k)^2$, for which integration by parts in $y$
and the estimates of Remark \ref{eboundsrmk} yield integrals of the
form
\begin{equation*}
\int_0^t \int_{-\infty}^0 (t-s)^{-1/2} e^{-\frac{(y+a_j^-
(t-s))^2}{M(t-s)}} (1+s)^{-1} e^{-\frac{(y - a_k^- s)^2}{Ms}} dyds,
\end{equation*}
with $a_j^- > 0$ (the excited terms correspond with mass convecting
into the shock layer) and $a_k < 0$ (diffusion waves care mass away
from the shock layer).  According to Lemma 6 of \cite{HZ}, we have
the relation
\begin{equation*}
\begin{aligned}
&e^{-\frac{(y+a_j^- (t-s))^2}{M(t-s)}} e^{-\frac{(y-a_k^- s)^2}{Ms}}
=e^{-\frac{(a_j^- (t-s) + a_k^- s)^2}{Mt}} \\
&\times e^{-\frac{t}{Ms(t-s)}(y+\frac{(a_j^- + a_k^-)(t-s)s}{t})^2},
\end{aligned}
\end{equation*}
from which integration over $y$ leads to an estimate of the form
\begin{equation*}
Ct^{-1/2} \int_0^t (1+s)^{-1/2} e^{-\frac{(a_j^- (t-s)+a_k^-
s)^2}{Mt}} ds.
\end{equation*}
In either the case $s \in [0,t/\gamma]$ or $s \in [t-t/\gamma, t]$,
for $\gamma$ sufficiently large, we have
\begin{equation*}
|a_j^- (t-s) + a_k^- s| \ge \eta t,
\end{equation*}
for some $\eta > 0$, through which we have exponential decay in $t$.
In the case $s \in [t/\gamma, t-t/\gamma]$, we integrate the kernel
in $s$ to obtain an estimate by $C (1+t)^{-1/2}$.

For the final estimate in Lemma \ref{nonlinearestimates}, we have
for $y<0$ integrals of the form
\begin{equation*}
\int_0^t \int_{-\infty}^0 (t-s)^{-1} e^{-\frac{(y+a_k^-
(t-s))^2}{M(t-s)}} \Psi (y,s) dy ds,
\end{equation*}
where $a_k^- < 0$.  In the case of nonlinearities $(1+s)^{-1/4}
s^{-1/2} (\bar\psi_1 + \psi_2 + \varphi)$, we can proceed by setting
$x = 0$ in the scattering estimates to obtain an estimate by
$(1+t)^{-5/4}$. In the case of nonlinearity $(1+s)^{-1} s^{-1/2}
(1+|y|)^{-1/2}$, we have the estimate
\begin{equation*}
\begin{aligned}
e^{-\frac{(y+a_k^- (t-s))^2}{M(t-s)}} (1+|y|)^{-1/2} \le C
e^{-\frac{(y+a_k^- (t-s))^2}{\bar{M}(t-s)}} (1+|y|+(t-s))^{-1/2},
\end{aligned}
\end{equation*}
from which we immediately obtain an estimate by
\begin{equation*}
\begin{aligned}
C_1& t^{-1/2} (1+t)^{-1/2} \int_0^{t/2} (1+s)^{-1} s^{-1/2} ds \\
&+ C_2 (1+t)^{-1} t^{-1/2} \int_{t/2}^t (t-s)^{-1/2} (1+(t-s))^{-1/2} ds \\
&\le C (1+t)^{-1}.
\end{aligned}
\end{equation*}
Finally, in the case of nonlinarity $(1+s)^{-1} e^{-\eta |y|}$, we
have the estimate
\begin{equation*}
e^{-\frac{(y+a_k^- (t-s))^2}{M(t-s)}} e^{-\eta |y|} \le C e^{-\eta_1
(t-s)} e^{-\eta_2 |y|},
\end{equation*}
from which we immediately obtain an estimate by
\begin{equation*}
\begin{aligned}
C_1 &e^{-\frac{\eta_1}{4}t} \int_0^{t/2} e^{-\frac{\eta_1}{2} (t-s)}
ds
+ C_2 (1+t)^{-1} \int_{t/2}^t e^{-\eta_1 (t-s)} ds \\
&\le C (1+t)^{-1}.
\end{aligned}
\end{equation*}

This concludes the proof of Lemma \ref{nonlinearestimates} \hfill
$\square$

\medskip
\noindent {\bf Proof of Lemma \ref{Hnonlinear}.} To show (\ref{Hf})
we need to estimate
$$\int_0^t\int_{-\infty}^{+\infty}\mathcal{R}_j^*(x) \mathcal{O}(e^{-\eta_0 (t-s)}) \delta_{x-\bar a_j^*
(t-s)}(-y) \mathcal{L}_j^{*t}(y)\Upsilon(y,s) dyds$$
$$\le C\int_0^t (e^{-\eta_0 (t-s)})|f(-x+\bar a_j^*
(t-s),s)|ds$$
 by the boundedness of $\mathcal{R}_j^*$ and  $\mathcal{L}_j^{*t}$. A typical term of
 $\bar\psi_1$
 is  a term of the form
$(1+t)^{-\frac12}(1+|x-a_i^-t|)^{-\frac34}$ (or with $a_i^-$
replaced by $a_i^+$). Hence for the terms coming from $\bar\psi_1$
the above is of the order
$$\int_0^t e^{-\eta_0(t-s)} s^{-\frac12}(1+s)^{-\frac12}(1+|x-\bar a_j^* (t-s)- a_i^-s |)^{-\frac34} ds $$
$$=\int_0^t e^{-\eta_0(t-s)} s^{-\frac12}(1+s)^{-\frac12}(1+|x- a_i^-t -(\bar a_j^*-a_i^-) (t-s)|)^{-\frac34} ds;$$
now, using (\ref{abineq}), this is smaller than
$$ \int_0^t e^{-\eta_0(t-s)} s^{-\frac12}(1+s)^{-\frac12}(1+|x- a_i^-t |)^{-\frac34}(|1+(\bar a_j^*-a_i^-) (t-s)|)^{\frac34} ds. $$
Notice that $(\bar a_j^*-a_i^-)\le C$. Now we use half of $\eta$ to
neutralize $t-s$, and to get
$$(1+|x- a_i^-t |)^{-\frac34} \int_0^t e^{-\frac{\eta_0}2(t-s)} s^{-\frac12}(1+s)^{-\frac12}ds$$
$$\le (1+|x- a_i^-t |)^{-\frac34}(1+t)^{-1},$$
 obviously absorbable in $\bar\psi_1+\psi_2+\alpha.$ For the terms coming from $\psi_2$, we notice
 that, by inequality
$a^2+b^2\ge 2ab,$
  $$(1+|y-a_i^-s|+\sqrt s)^{-\frac 32}\le C(1+s)^{-\frac38}(1+|y-a_i^-s|)^{-\frac34}.$$  Hence applying the same
  calculations as for $\bar\psi_1$ gives us  an estimate of $(1+|x- a_i^-t |)^{-\frac34}(1+t)^{-\frac78}.$
  For  $\alpha$, we apply the
 same procedure, and using again (\ref{abineq}), we have an estimate by
 \begin{equation*}
 \begin{aligned}
C\int_0^t &e^{-\eta_0(t-s)} s^{-\frac12}(1+s)^{-\frac34}(1+|x-\bar
a_j^*
(t-s)|)^{-\frac12} ds \\
&\le C\int_0^t e^{-\eta_0(t-s)}
s^{-\frac12}(1+s)^{-\frac34}(1+|x|)^{-\frac12}
(1+|\bar a_j^* (t-s)|)^{\frac12}ds \\
&\le C(1+|x|)^{-\frac12}\int_0^t e^{-\frac{\eta_0}2(t-s)}
(1+s)^{-\frac54} ds \\
&\le C(1+|x|)^{-\frac12} (1+t)^{-1}.
\end{aligned}
\end{equation*}
  This finishes the proof of
(\ref{Hf}) and (\ref{Hxf}). The proof of (\ref{Hphi}) is similar,
since all the terms in $\Phi$ are smaller than  terms looking like
 $C(1+t)^{-\frac32}e^{\frac{-(x-a_i^- t)^2}{Mt}}$, which
is bounded by $C(1+t)^{-\frac34}(1+ |x-a_i^-t|+\sqrt t)^{-\frac
32}.$

This finishes the proof of Lemma \ref{Hnonlinear}.
  \hfill $\square$

\end{document}